\documentclass[11pt]{article}
\usepackage{times,amsmath,amsfonts,amsthm,amssymb,eucal}
\usepackage{algorithm}
\usepackage{algorithmic}
\usepackage{graphicx,ifthen}
\usepackage{wrapfig}
\usepackage{color}
\newtheorem{theorem}{Theorem}[section]
\newtheorem{lemma}{Lemma}[section]

\newtheorem{remark}{Remark}[section]
\newtheorem{corollary}{Corollary}[section]
\newtheorem{proposition}{Proposition}[section]

\def\bs{\boldsymbol}

\def\ol{\overline}
\begin{document}

%\begin{frontmatter}
%\title{Nonparametric estimation of dynamics of monotone trajectories}
%\runtitle{Dynamics of monotone trajectories}

\begin{center}

\Large{\bf Nonparametric estimation of dynamics of monotone trajectories}

\vskip.1in

\large{Debashis Paul, Jie Peng and Prabir Burman}\footnote{Paul's  research is partially supported by the NSF grants DMR-10-35468 and
DMS-11-06690. Peng's research is partially supported by the NSF grant DMS-10-01256.
Burman's research is partially supported by the NSF grant DMS-09-07622.}

\vskip.1in

\textit{\small Department of Statistics, University of California, Davis}

\end{center}

%\begin{aug}
%\author{\fnms{Debashis} \snm{Paul}\thanksref{t1}\ead[label=e1]{debpaul@ucdavis.edu},
%\fnms{Jie} \snm{Peng}\thanksref{t2}\ead[label=e2]{jiepeng@ucdavis.edu},  and \fnms{Prabir} %\snm{Burman}\thanksref{t3}\ead[label=e3]{pburman@ucdavis.edu}}

%\thankstext{t1}{Paul's  research is partially supported by the NSF grants DMR-10-35468 and
%DMS-11-06690.}
%\thankstext{t2}{Peng's research is partially supported by the NSF grant DMS-10-01256.}

%\thankstext{t3}{Burman's research is partially supported by the NSF grant DMS-09-07622.}

%\runauthor{Paul, Peng and Burman}

%\affiliation{University of California, Davis}
%\address{Department of Statistics\\
%University of California\\
%Davis, CA 95656, U. S. A.\\
%\printead{e1}\\
%\phantom{E-mail:\ }\printead*{e2}\\
%\phantom{E-mail:\ }\printead*{e3}}
%\end{aug}

\begin{abstract}
We study a class of nonlinear nonparametric inverse problems.
Specifically, we propose a nonparametric estimator of the dynamics of a monotonically
increasing trajectory defined on a finite time interval. Under suitable regularity conditions,
we prove consistency of the proposed estimator and show that in terms of $L^2$-loss,
the optimal rate of convergence for the proposed estimator is the same
as that for the estimation of the derivative of a trajectory.
This is a new contribution to the area of nonlinear nonparametric inverse problems.
We conduct a simulation study to examine the finite sample behavior
of the proposed estimator and apply it to the Berkeley growth data.
\end{abstract}

%\begin{keyword}[class=AMS]
%\kwd[Primary ]{62G08}
%\kwd[; secondary ]{62G20}
%\end{keyword}

%\begin{keyword}
%\kwd{autonomous differential equation}
%\kwd{nonlinear inverse problem}
%\kwd{monotone trajectory}
%\kwd{nonparametric estimation}
%\kwd{perturbation theory}
%\kwd{spline}
%\end{keyword}

%\end{frontmatter}


\section{Introduction}

Monotone trajectories describing the evolution of state(s) over time appear
widely in scientific studies. The most widely studied are probably growth of
organisms such as humans or plants (Milani, 2000; Erickson, 1976; Silk and
Erickson, 1979). There are many parametric models for describing the features
of growth curves, particularly in human growth (Hauspie et al., 1980; Milani,
2000). Most of these works focus on modeling the trajectories themselves or
modeling the rate of change, i.e., the derivative of the trajectories.
Other examples of monotone trajectories appear in
population dynamics under negligible resource constraints (Turchin, 2003), in
dose-response analysis in pharmacokinetics (Kelly and Rice, 1990), in auction
price dynamics in e-Commerce (Jank and Shmueli, 2006; Wang et al., 2008; Liu
and M\"{u}ller, 2009), and in analysis of trajectories of aircrafts after
take-off (Nicol, 2013). Some of these works  are looking at function estimation
with monotonic constraints and some of them are taking a functional data
analysis approach.

In contrast, our goal here is to estimate the functional relationship between
the rate of change and the state, i.e., the dynamics of the trajectory, through
a nonparametric model.
%Relating the rate of change with the state provides
%important insights about the underlying system and is often of great scientific
%interests.
Many systems such as growth of organisms or economic activity of a
country/region are intrinsically dynamic in nature (cf. Ljung and Glad, 1994).
A dynamics model provides a mechanistic description of  the system rather than
a purely phenomenological one. Moreover, due to insufficient scientific
knowledge, quite often there is a need for nonparametric modeling of the
dynamical system. In addition, nonparametric fits can be used to develop
measures of goodness-of-fit for hypothesized parametric models.
%However, there is little work on this aspect in the literature.
%Specifically,  we propose a nonparametric estimation of  an autonomous ordinary differential equation

There is a large literature in modeling continuous time smooth dynamical
systems through systems of parametric differential equations (see, e.g.,
Perthame, 2007, Strogatz, 2001). These methods have been used to model HIV
dynamics (Wu, Ding and DeGruttola, 1998; Wu and Ding, 1999; Xia, 2003; Chen and
Wu, 2008a, 2008b),  the dynamic behavior of gene regulation networks (Gardner
\textit{et al.}, 2003; Cao and Zhao, 2008), etc. Approaches for fitting a
parametric dynamics model include the maximum likelihood or nonlinear least
squares.
%approach in which the solution of the ODE needs to be evaluated
%given the parameter at all the time points and this needs to be done numerically
%except for a linear ODE.
A recent  approach proposed by Ramsay \textit{et al.} (2007) and Cao \textit{et
al.} (2008) for parametric ordinary differential equations is based on the idea
of balancing the model fit and the goodness of fit of the trajectories
simultaneously.

Another popular approach to fit dynamics models is a two-stage procedure (Chen
and Wu, 2008a, 2008b;, Varah, 1982), where the trajectories and their
derivatives are first estimated nonparametrically  and then the dynamics is
fitted by regressing the fitted derivatives to the  fitted trajectories.
%Variants of this procedure have been studied by Chen and Wu (2008a, 2008b), who estimated differential equations
%with known functional forms and nonparametric time-dependent coefficients
%through. This method has precedence in Varah (1982)
%who used cubic splines for fitting the sample trajectories in the first
%stage of the fitting procedure.
The two-stage approach can be easily adapted to estimate a nonparametric
dynamics model. However, their performance is unsatisfactory due the difficulty
of resolving the bias-variance trade-off in a data dependent way. Brunel (2008)
gives a comprehensive theoretical analysis of such an approach.
Very recently, Hall and Ma (2014) proposed a one-step estimation procedure that
mitigates some of the inefficiencies of two-stage estimators. However, this approach
does not seem to extend naturally to estimate nonparametric dynamical systems.

There is also an extensive literature on the nonparametric estimation of monotone
functions, e.g.,  Brunk (1970), Wright and Wegman (1980),  Mammen (1991),
Ramsay (1988, 1998). However, most methods in this field are not concerned with
the estimation of the gradient function, except for Ramsay (1998) where the
unknown function is modeled in terms of a second order differential equation
and  a smoothed estimate of its gradient is obtained as a byproduct.

%In this paper, we propose a nonparametric
%estimator of the gradient function of an autonomous nonlinear dynamical system.
%based on noisy observations from a monotonically increasing trajectory on a finite time interval.
%We  study its theoretical properties and numerical behavior. \textcolor{red}{As far as we know, this is the first work in literature %to .....}
%Most of the existing approaches assume known
%functional forms of the gradient function,
%i.e., the function relating $X'(t)$ to $X(t)$ and $t$.
%and many of them are developed to fit the differential equations for a single subject.
%A notable exception is the works of Chen and Wu (2008a, 2008b) and Xue, Miao and Wu (2010)
%who allow the parameters describing the ODE to be smooth, time-dependent functions.

A key observation of estimating the dynamics of monotone trajectories
is that for any smooth monotone trajectory, its
dynamics can be described by a first order autonomous differential equation. Specifically,
if $X(t)$ is positive, strictly monotone and differentiable on a finite time
interval, then we can express
\begin{equation}\label{eq:autonomous_representation}
X'(t) = (X'~o~ X^{-1})(X(t)) = g(X(t)), ~~~t \in [0,1]
\end{equation}
where $g = X' ~o~ X^{-1}$ is the gradient function.
In this paper, we estimate the unknown gradient function $g$
nonparametrically from discrete noisy observations of $X$. Specifically, we model the
gradient function by a basis representation where the number of basis functions
grow with the sample size. We adopt a nonlinear least squares framework for model fitting.
%in which the solution of the ODE
%is evaluated explicitly given the coefficients of the basis representation.
%The
%proposed method is applicable to estimate the dynamics of any one-dimensional variable
%where the true trajectory is a positive, monotone function, e.g. growth curves.
%We also present a method for selecting the number of basis functions by an
%approximate leave-one-curve out cross-validation method.
We then carry out a detailed theoretical analysis and derive the rate of
convergence of the proposed estimator.
%We present simulation studies and apply the method to Berkeley Growth data (Tuddenham and Snyder, 1954).

%At this point, we would like to clarify the distinction of our
%methodology from the extensive literature on the nonparametric estimation of
%monotone functions which falls in the category of shape constrained regression.
%Important works in this direction include Brunk (1970), Banerjee and Wellner
%(2001), Kelly and Rice (1990), Hall and Huang (2001), Mammen (1991), Mammen and
%Thomas-Agnan (1999), Ramsay (1988, 1998), Mukerjee (1988) and Wright and Wegman
%(1980). Bouzas and Ruiz-Fuentes (2011) considered population level modeling of
%monotone trajectories based on a functional data analysis. Most of these
%methods are not concerned with the estimation of the gradient function of the
%monotone trajectory, Ramsay (1998) being a notable exception in that there the
%unknown regression function is modeled in terms of a second order differential
%equation and consequently a smoothed version of the gradient function is
%obtained as a byproduct. In contrast, our objective is to estimate the gradient
%function of the trajectory itself, which is a nonparametric nonlinear inverse
%problem as opposed to the monotonic regression problem which has a different
%characteristic. Of course, our gradient estimate can be used to estimate
%monotonic fits to the trajectories, even though no further theoretical analysis
%is carried out for the latter object.

We  now highlight the major contributions of this work.
Although there is a large literature on linear nonparametric inverse problems
(Cavalier \textit{et al.}, 2004; Cavalier, 2008; Donoho, 1995; Johnstone
\textit{et al.}, 2004), especially on the nonparametric estimation of the
derivative of a curve (Gasser and M\"{u}ller, 1984; M\"{u}ller et al., 1987;
Fan and Gijbels, 1996),
%Also, there is a substantial literature on nonlinear parametric inverse problems.
there is little theoretical development on nonlinear nonparametric inverse
problems. Thus, our work makes a new contribution to this important area.
In this paper, we first quantify the degree of ill-posedness  of the estimation
of the gradient function $g$
%by finding the order of the
%condition number of the equivalent of the information matrix
as the number of basis functions grow to infinity. We then use this result to show that if $g$ is
$p$ times differentiable then the $L^2$-risk of the proposed estimator has the same optimal rate of convergence,
viz., $O(n^{-2p/(2p+3)})$, as that of the estimator of the derivative of a trajectory
assuming that the latter is $p+1$ times differentiable. In Section \ref{sec:discussion},
we show that the optimal rate of the proposed estimator is indeed the minimax rate for estimation
of $g$ under $L^2$ loss if the class of estimators is restricted to be uniformly Lipschitz.
%We conjecture that this rate optimality holds even if we drop the Lipschitz requirement on the estimators.
In the rest of the paper, unless otherwise specified, the phrase ``optimal rate'' refers to the best rate
of convergence of the proposed estimator.

Among the few instances of nonparametric modeling
of the gradient function known to us, Xue, Miao and Wu (2010) dealt with a related but different
problem of estimating a parametric ODE with time-varying parameters, where the latter
are modeled as unknown smooth functions of time. In a work most
closely related to ours, Wu et al. (2014) proposed a sparse additive
model for describing the dynamics of a multivariate state vector and developed a combination
of two-stage smoothing and sparse penalization for fitting the model. Their model can be
seen as a multi-dimensional generalization of the autonomous ODE model studied here.
In their paper, while deriving the risk bounds, it is assumed that whenever the gradient
function $g$ is $p$ times differentiable, the state $X$ is at least $3p+1$ times differentiable.
However, due to the representation $g = X' ~o~ X^{-1}$, it follows that $g$ is  $p$
times differentiable if and only if $X$ is $p+1$ times differentiable.
Therefore, at least for the one-dimensional state variable case, the assumptions made
in Wu et al. (2014) are not satisfied in reality if $p$ indeed denotes the maximal
order of smoothness of $g$. This indicates that the rate of convergence
their estimator of $g$ is not optimal for the current problem.
It is also instructive to note that, due to the assumption about the additional degree
of smoothness of the state variable, Wu et al. (2014)
did not encounter the technical challenge posed by the ill-posedness of the problem.

%In particular, we show that if the gradient function is sufficiently smooth,
%then the $L^2$-risk of the proposed estimator has the same asymptotic convergence
%rate as that of the derivative estimator in usual nonparametric regression problems,
%even though the latter is a linear inverse problem, and consequently, easier
%to analyze theoretically.

%For more details, see equation (\ref{eq:G_star_def}).
%We also formulate a stable computational procedure for obtaining the
%estimator.

The rest of the paper is organized as follows. In Section \ref{sec:model}, we
briefly describe the model and the estimation procedure. We
present the main theoretical results in Section \ref{sec:theory} and outline
the main steps of the proof in Section \ref{sec:proofs}. We present a
simulation study in Section \ref{sec:simulation} and an application to the
Berkeley growth data in Section \ref{sec:real}. We discuss the optimality
of the estimation of $g$ in Section \ref{sec:discussion}. Some proof details
are provided in the Appendix (Section \ref{sec:Appendix}). Some derivations
and graphical summaries are provided in the Supplementary Material (SM).

\section{Model}\label{sec:model}

The class of models studied in this paper is of the form:
\begin{equation}\label{eq:basic}
X'(t) = g(X(t)), ~~ X(0)=x_0, ~~~ t\in [0,1],
\end{equation}
where $g$ is an unknown smooth function which is assumed to be positive on the
range of $\{X(t):t \in [0,1]\}$. Therefore, the sample trajectory $X(t)$ is a
strictly increasing function of time $t$. The observations are
\begin{equation}\label{eq:data_model}
Y_j = X(t_j) + \varepsilon_j,~~j=1,\ldots,n,
\end{equation}
where $0 \leq t_1,\ldots, t_n \leq 1$ are
observation times. The noise terms $\varepsilon_j$'s are assumed to be i.i.d. with mean 0 and
variance $\sigma_{\varepsilon}^2 > 0$.

Our goal is to estimate the gradient function $g$ based on the
observed data $Y_j$s. We propose to  approximate $g$ through a basis
representation:
\begin{equation}\label{eq:basis}
g(x) \approx g_{\bs\beta} := \sum_{k=1}^M \beta_k \phi_{k,M}(x),
\end{equation}
where $\{\phi_{k,M}(\cdot)\}_{k=1}^M$ is a set of linearly independent
compactly supported smooth functions. Henceforth, we use
$\phi_k$ to denote $\phi_{k,M}$.

%The
%theoretical derivations in this paper are based on modeling $g$ in terms of
%cubic B-spline basis functions with equally spaced knots, even though other
%choices of compactly supported basis functions with sufficient degrees of
%smoothness is also possible.

%\textcolor{blue}{This formulation is
%related to the formulation of semiparametric autonomous dynamics model for
%plant root growth studied by Paul, Peng and Burman (2011). However, no
%theoretical analysis of the procedure was presented in that paper. The
%theoretical derivations in this paper are based on modeling $g$ in terms of
%cubic B-spline basis functions with equally spaced knots, even though other
%choices of compactly supported basis functions with sufficient degrees of
%smoothness is also possible.}
%If the trajectory $X$ is increasing (without loss
%of generality) then $g$ is positive on its domain. However, this positivity
%constraint is not imposed in our model. A simple approach that would allow such
%a formulation would be to represent $\log g$ rather than $g$ in a spline basis.
%But our present formulation is much more convenient both in terms of numerical
%implementation and theoretical analysis since the perturbation analysis of the
%corresponding ODE is considerably simpler. Moreover, the numerical algorithm is
%very stable and the consistency results imply that the fitted $g$ is positive
%with probability tending to one if the true trajectory is strictly increasing.}

We now describe the estimation procedure. For the time being, assume that we
observe the two endpoints $x_0 = X(0)$ and $x_1 = X(1)$ noiselessly and so the
combined support of $\{\phi_{1},\ldots,\phi_{M}\}$ is the interval $[x_0,x_1]$.
%In practice, we take the combined support of the basis functions to be slightly bigger interval
%containing $[x_0,x_1]$, for reasons to be explained later.
Given any $\bs\beta := (\beta_1,\ldots,\beta_M)$ so that $g_{\bs\beta}$ is
positive on the support of $\{\phi_{k}(\cdot)\}_{k=1}^M$, we can solve the
initial value problem
\begin{equation}\label{eq:initial_general}
x'(t) = g_{\bs\beta}(x(t)), ~~~t \in [0,1], ~~~x(0) = x_0
\end{equation}
to obtain the corresponding trajectory $X(t;\bs\beta)$.
Define the $L^2$ loss function:
\begin{equation}\label{eq:log-like}
L(\bs\beta):= \sum_{j=1}^{n} (Y_j - X(t_j;\bs\beta))^2.
\end{equation}
Then the proposed estimator of $g$ is defined as
\begin{equation}\label{eq:hat_g_beta}
\widehat{g}(x) := g_{\widehat{\bs\beta}}(x) = \sum_{k=1}^M \widehat{\beta}_k \phi_{k,M}(x),
~~~{\rm where} ~~ \widehat{\bs{\beta}} := {\rm arg}\min_{\bs\beta \in \mathbb{R}^M}L(\bs\beta).
\end{equation}
Minimization of $L(\bs\beta)$ is  a nonlinear least squares problem. We propose to
use a Levenberg-Marquardt iterative updating scheme. Since this requires
evaluating the trajectory $X(t;\bs\beta)$ and its derivative with respect to
$\bs\beta$, given the current estimate of $\bs\beta$, we solve the
corresponding differential equations numerically by using the 4-th order
Runge-Kutta method. More details are given in the Appendix.
%Levenberg-Marquardt scheme also requires specifying an initial estimate for
%$\bs\beta$.  This can be obtained by separately estimating $X(\cdot)$ and $X'(\cdot)$ by,
%e.g., local polynomial regression, and then obtaining linear regression estimate of $\bs\beta$
%by treating $\widehat X'$ as the response and
%$\phi_{1,M}(\widehat X),\ldots,\phi_{M,M}(\widehat X)$ as predictors, as is done in Section
%\ref{subsec:two_stage_regerssion}.
%This can be obtained by a nonparametric regression estimate of $g$
%with the empirical derivatives $(X(t_{j+1})-X(t_{j}))/(t_{j+1}-t_j)$'s as the
%response and empirical observations $(X(t_{j+1})+X(t_j))/2$'s as the predictor
%and then projecting the estimate onto the basis $\{\phi_{1},\ldots,\phi_{M}\}$.
Finally, the number of basis $M$ is selected through an approximate
cross-validation score. A fitting procedure using similar
techniques is studied in Paul et al. (2011) in a different context.

In practice,  the initial value $x_0 = X(0)$ and  the right boundary $x_1 =
X(1)$ may not be observed or may be observed with noise. The choice of the
endpoints of the combined support of the basis functions then becomes a
delicate matter. This is because evaluation of the trajectory is an initial
value problem, so error in $x_0$ propagates throughout the time domain.
%Moreover, uncertainty in $x_1$ leads to boundary effects which is of a
%different nature than the boundary issues in ordinary non-parametric regression
%problems.
We discuss this in more details in the Appendix (particularly, see
Figure \ref{fig:trajectory_envelop}).

In the following, we  propose a modified estimation procedure when $x_0$ and
$x_1$ are unknown. The basic idea is to first estimate the trajectory at the
endpoints of a slightly smaller time interval $[\delta, 1-\delta]$ for a small
positive constant $\delta$, and then estimate the gradient function using data
falling within this time interval. Throughout the paper, $\delta$ is treated as
a fixed quantity. In practice, we may select $\delta$ to be the time point such
that about 5\% of the data fall in the intervals $[0,\delta]$ and
$[1-\delta,1]$. Too small a value of $\delta$ may cause distortions of the
estimated $g$ at the boundaries.

%Specifically, instead of estimating $g$ on the interval $[x_0,x_1]$, we
%estimate it on $[x_{0,\delta}, x_{1,\delta}]$ where $x_{0,\delta} := X(\delta)$
%and $x_{1,\delta} = X(1-\delta)$ and $\delta$ is a small positive number.

We first  obtain  nonparametric estimates of $x_{0,\delta} := X(\delta)$ and
$x_{1,\delta}:= X(1-\delta)$, denoted by $\widehat{x}_0$ and $\widehat{x}_1$,
respectively. We then define
%\begin{equation}\label{eq:eta}
%\end{equation}
$x_{0,M} = \widehat {x}_0 -\eta_M$ and $x_{1,M} = \widehat {x}_1 + \eta_M$,
where $\eta_M$ is a small positive number satisfying $\eta_M = o(M^{-1})$
which implies that $x_0 < x_{0,M} < x_{1,M} < x_1$ as $n$ goes to infinity. At
the same time, $\eta_M$ should be large enough so that
$\max_{j=0,1}|x_{j,\delta} - \widehat x_j| = o_P(\eta_M)$  which ensures that
$x_{0,M} < x_{0,\delta} < x_{1,\delta} < x_{1,M}$ and
$\max_{j=0,1}|x_{j,\delta}-x_{j,M}| = O_P(\eta_M) = o_P(M^{-1})$ as $n$ goes to infinity. For some
technical considerations, to be utilized later, we also want $\eta_M \gg M^{-3/2}$. In practice, we may
select $\eta_M$ to be $\min\{M^{-3/2}\log n,s_M/\log n \}$ where $s_M$ is the length of the
smallest support among the basis functions $\{\phi_1,\ldots,\phi_M\}$. For more
details on how to obtain $\widehat x_j$, $j=0,1$, see Lemma
\ref{lem:consistency_x_j_hat} in Section \ref{sec:theory}.  In addition, we
also assume that $\widehat x_j$, $j=0,1$ are estimated from a sample
independent from that used in estimating $\bs\beta$. This can be easily
achieved in practice by sub-sampling of the measurements. This assumption
enables us to prove the consistency result (in Section \ref{sec:theory})
conditionally on $\widehat x_j$, $j=0,1$ and treating them as nonrandom
sequences converging to $x_{j,\delta}$, $j=0,1$.

We  then set the combined support of the basis functions
$\{\phi_{k,M}\}_{k=1}^M$ as the interval $[x_{0,M}, x_{1,M}]$, and use the
following modified  loss function  to derive an estimator for $g$:
\begin{equation}\label{eq:loss_modified}
\tilde L_\delta(\bs\beta) = \sum_{j=1}^n (Y_j - X(t_j;\bs\beta,\widehat x_0))^2 \mathbf{1}_{[\delta,1-\delta]}(t_j),
\end{equation}
where $X(t;\bs\beta,a)$ denotes the integral curve of the ODE
\begin{equation}\label{eq:initial_beta_delta}
x'(t) = g_{\bs\beta}(x(t)), ~~~t \in [\delta,1-\delta], ~~~x(\delta) = a.
\end{equation}
The estimated $\widehat g$ is through minimizing the above loss function with
respect to $\bs\beta$ (equation (\ref{eq:hat_g_beta}) with $L$ replaced by $\tilde L_\delta$).

\section{Consistency}\label{sec:theory}

In this section, we discuss the consistency of the estimator $\widehat{g}$
defined by the loss function (\ref{eq:loss_modified}). The asymptotic framework is that the number of
basis functions $M$ goes to infinity together with the number of measurements
$n$.  The consistency of the estimator $\widehat g$ over
$[x_{0,\delta},x_{1,\delta}]$ is formulated in terms of the $L^2$-loss as:
$$
\int_{x_{0,\delta}}^{x_{1,\delta}} |\widehat g(u) - g(u)|^2 du \longrightarrow 0,
~~~ \hbox{in probability as $n \rightarrow \infty$}.
$$
In Theorem \ref{thm:optimal_rate} we derive a bound on the rate of convergence of $\widehat g$ in terms of the $L^2$-loss
as $n,M \to \infty$ that depends upon the degree of smoothness of $g$. Specifically, the optimal rate is
$O_P(n^{-2p/(2p+3)})$ for $p \geq 4$.

\subsection{Assumptions}

The following assumptions are made on the model.
\begin{itemize}
\item[{\bf A1}]
$g \in C^p(D)$, and $g > 0$ on $D$ for some
integer $p \geq 3$, where $D$ is an open interval containing $[x_0,x_1]$.
\item[{\bf A2}]
The collection of basis functions $\Phi_M :=
\{\phi_{1,M},\ldots, \phi_{M,M}\}$ satisfies:
\begin{itemize}
\item[(i)]
$\phi_{k,M}$'s have unit $L^2$ norm;
\item[(ii)]
the combined support of $\Phi_M$ is $D_0 \equiv D_{0,M} := [x_{0,M},x_{1,M}]$ and
for every $k$, the length of the support of $\phi_{k,M}$ is
$O(M^{-1})$;
\item[(iii)]
$\phi_{k,M} \in C^2(D_0)$ for all $k$;
\item[(iv)]
$\sup_{x \in D_0} \sum_{k=1}^M |\phi_{k,M}^{(j)}(x)|^2
= O(M^{1+2j})$, for $j=0,1,2$;
\item[(v)]
the Gram matrix $\mathbf{G}_{\Phi_M} := ((\int_{x_{0,M}}^{x_{1,M}} \phi_{k,M}(u) \phi_{l,M}(u)du))_{k,l=1}^M$
is such that there exist constants $0 < \underline{c} \leq \overline{c} < \infty$, not depending on $M$
such that $\underline{c} \leq \lambda_{min}(\mathbf{G}_{\Phi_M}) \leq \lambda_{max}(\mathbf{G}_{\Phi_M}) \leq \overline{c}$ for all $M$;
\item[(vi)]
for every $M$, there is a $\bs{\beta}^*\in \mathbb{R}^M$ such that
$\sup_{t \in [\delta,1-\delta]}|X_g(t) - X(t;\bs\beta^*)| = O(M^{-(p+1)})$ and
$\sup_{u\in [x_{0,\delta},x_{1,\delta}]} |g^{(j)}(u) - g_{\bs\beta^*}^{(j)}(u)| = O(M^{-p+j})$
for $j=0,1,2$, where $g_{\bs\beta} = \sum_{k=1}^M \beta_k
\phi_{k,M}$ and $X(t;\bs\beta) \equiv X(t;\bs\beta,x_{0,\delta})$ with $X(t;\bs\beta,a)$ as in
(\ref{eq:initial_beta_delta}).
\end{itemize}
\item[{\bf A3}]
Time points $\{t_j\}_{j=1}^{n}$ are realizations of
$\{T_j\}_{j=1}^{n}$, where $T_j$'s  are i.i.d.
from a continuous distribution $F_T$ supported on $[0,1]$ with a
density $f_T$ satisfying $\underline{c}' \leq f_T \leq \overline{c}' $ for some $0 < \underline{c}'
\leq \overline{c}' < \infty$.
\item[{\bf A4}]
The noise $\varepsilon_j$'s are i.i.d. sub-Gaussian random variables (cf. Vershynin, 2010)
with mean 0 and variance $\sigma_\varepsilon^2 > 0$.
\end{itemize}
We give brief explanations of these assumptions. {\bf A1} ensures
sufficient smoothness of the solution paths of the differential equation
(\ref{eq:basic}). Also by {\bf A1}, $X_g(\cdot)$ is $p+1$ times continuously differentiable on $D$.
Assumptions (i) to (v) of {\bf A2} are satisfied by  B-spline
basis, rescaled to have unit norm, of order $\geq 3$, with equally spaced knots .
Define
\begin{equation}\label{eq:xi_n}
\xi_n := \sqrt{\log n} n^{-\frac{p+1}{2p+3}},
\end{equation}
which is used in determining the rates of convergence of the estimator.
Thus, by making use of {\bf A4}, we get the following results  with respect to the estimates of $x_{0,\delta}$ and $x_{1,\delta})$
(cf. Fan and Gijbels, 1996).
\begin{lemma}\label{lem:consistency_x_j_hat}
Suppose that {\bf A1} and {\bf A4} hold. Consider using a kernel of sufficient degree of
smoothness to obtain estimates $\widehat x_j$ for $x_{j,\delta}$, $j=1,2$,
through local polynomial method with bandwidth of order $n^{-1/(2p+3)}$.
Define $d_n := \max_{j=0,1}|\widehat x_j - x_{j,\delta}|$. Then $d_n = O_P(n^{-(p+1)/(2p+3)})$
and given $\eta > 0$, there exists $C(\eta) > 0$
such that $d_n \leq C(\eta) \xi_n$ with probability at least $1-n^{-\eta}$, where $\xi_n$ is as in
(\ref{eq:xi_n}).
\end{lemma}
If  $M=O((n/\log n)^{1/7})$ (as in Theorem
\ref{thm:consistency}) and $M^{-3/2}\ll \eta_M \ll M^{-1}$ for some $C > 0$, then we have $\xi_n
= o(\eta_M)$  as $n \to \infty$. This ensures that $D_0 = [x_{0,M},x_{1,M}]$
is within the interval $[x_0,x_1]$ a.s. for large enough $n$ and hence the properties of the function $g$ hold on $D_0$.
In addition,  $D_0$ contains the interval $[x_{0,\delta},x_{1,\delta}]$. Therefore condition (ii) in {\bf A2}
ensures that the combined support of the basis functions covers the range of
the data used in estimating $g$.

Conditions (i) to (v) of {\bf A2} are satisfied by many classes of basis functions,
including normalized B-spline basis of order $\geq 3$ with equally spaced knots in the interval
$[x_{0,M},x_{1,M}]$.  We show in Appendix B
%\ref{subsec:spline_approx}
that if the B-splines basis of order $\geq \max\{3,p-1\}$ with equally spaced knots in
$[x_{0,M},x_{1,M}]$, then (vi) of {\bf A2} is also satisfied.
Condition (vi) of {\bf A2} ensures that a solution $X(t;\bs\beta)$ of
(\ref{eq:initial_general}) on $t \in [\delta,1-\delta]$ exists for all
$\bs{\beta}$ sufficiently close to $\bs{\beta}^*$. This allows us to apply the
perturbation theory of differential equations to bound the fluctuations of the
sample paths when we perturb the parameter $\bs\beta$.

Assumption {\bf A3} on the randomness of the sample points allows us to work
with the random variables $\tilde T_j$ defined as $T_j$ conditional on $T_j \in
[\delta, 1-\delta]$ with conditional density $\tilde f_T$ given by $\tilde
f_T(t) = f_T(t)/(F_T(1-\delta) - F_T(\delta))$. The properties of $f_T$ ensure
that $\tilde f_T$ satisfies the same property on $[\delta,1-\delta]$ with
possibly modified values of the constants $c_1$ and $c_2$. It should be noted
that the key derivations leading to the consistency of $\widehat g$ are conditional
on $\mathbf{T}$ and therefore ${\bf A3}$ is only a convenient assumption for
describing the regularity of the time points. The asymptotic results (Theorems
\ref{thm:consistency} and \ref{thm:optimal_rate}) hold if instead of being randomly
distributed, the time points form a fixed regular grid, say, with equal spacing.

\subsection{Rate of convergence}

As mentioned earlier, the estimation of $g(\cdot)$ is a  nonlinear inverse
problem  since $X'(t)$ is not directly observable. In addition, this is also an
ill-posed estimation problem. Let
$X^{\bs{\beta}}(\cdot;\bs{\beta})$ be the partial derivative of
$X(\cdot;\bs\beta)$ with respect to $\bs{\beta}$, where $X(\cdot;\bs\beta)
\equiv X(\cdot;\bs\beta,x_{0,\delta})$ is the solution of
(\ref{eq:initial_beta_delta}) with $x(0) = x_{0,\delta}$. Let $\bs{\beta}^* \in
\mathbb{R}^M$ be as in {\bf A2}. Define
\begin{equation}\label{eq:G_star_def}
G_{*} := \mathbb{E}
\left(X^{\bs{\beta}}(\tilde T_1;\bs{\beta}^*)
(X^{\bs{\beta}}(\tilde T_1;\bs{\beta}^*))^T\right),
\end{equation}
where the expectation  is with  respect to the distribution of $\tilde T_1$.
Clearly $G_*$ is a positive semi-definite matrix. It becomes clear
from the analysis carried out later that the degree of ill-posedness of the
estimation problem is determined by the size of the operator norm of the matrix $G_{*}$
as a function of $M$. The following proposition gives a precise quantification of
the degree of ill-posedness. The situation here is in contrast with standard nonparametric function
estimation problems where the corresponding matrix is well-conditioned.

%We will show that there exists a sequence of positive numbers
%$\kappa_M$ (depending on $M$) such that, for every $M \geq M_0$ for some $M_0$
%depending on $p$,
%\begin{equation}\label{eq:G_star_condition}
%\parallel G_{*}^{-1}
%\parallel \leq \kappa_M.
%\end{equation}
%$\kappa_M$ increases to $\infty$
%with $M$ and its order as a function of $M$ quantifies the
%degree of ill-posedness of the estimation problem.
%Equation
%(\ref{eq:G_star_condition}) essentially quantifies the degree of
%ill-conditioning of the information matrix for $\bs{\beta}$.

\begin{proposition}\label{prop:kappa}
Assume that assumptions {\bf A1} to {\bf A3} hold with $p \geq 3$. Assume further that
(a) $\max_{j=0,1}|x_{j,M} - x_{j,\delta}| = o(M^{-1})$ (a.s.) and (b) $\min\{x_{1,M} - x_{1,\delta},x_{0,\delta}-x_{0,M}\}
\gg M^{-3/2}$. Then (a.s.)
\begin{equation}\label{eq:G_star_condition}
\parallel G_{*}^{-1} \parallel = O(M^2).
\end{equation}
\end{proposition}
By Lemma \ref{lem:consistency_x_j_hat} and the discussion that follows, under the condition
of Theorem \ref{thm:consistency}, (a) and (b) of Proposition \ref{prop:kappa} hold.

%\textcolor{blue}{Note that the condition $\max_{j=0,1}|x_{j,M} - x_{j,\delta}|
%= o(M^{-1})$ is ensured under the set up of Corollary
%\ref{cor:consistency_Bspline} stated below due to Lemma
%\ref{lem:consistency_x_j_hat}.}

%Moreover, a simple argument based on comparison
%of quadratic forms in the matrix $G_*$, by choosing an normalized version of
%the vector with coordinates alternating between $+1$ and $-1$, can be used to
%obtain a lower bound on $\kappa_M$ which is also of the order $O(M^2)$.
%Combining  this with Proposition \ref{prop:kappa}, we obtain that $\kappa_M
%\asymp M^2$.

We now state the main result on the consistency of the
estimate $\widehat g$.

%Note that, we have the result for a local minimizer
%of the objective function $\tilde L_\delta(\bs\beta)$ defined in (\ref{eq:loss_modified}).
%We also need to use a nonlinear optimization routine (Levenberg-Marquardt scheme in our
%case) to solve the optimization problem. We typically need a good initial estimate
%for achieving this. For this purpose, we use a two stage estimator discussed in
%Section \ref{subsec:two_stage_regerssion}.

%In the statement, we make the role of the
%ill-conditioning factor $\kappa_M$ explicit both in specifying the conditions
%on $M$ and in the rate of convergence. Then, via Proposition \ref{prop:kappa}
%we obtain the precise rate bounds  (Corollary \ref{cor:consistency_Bspline})
%when we use a normalized B-spline basis for representing $g$.

\begin{theorem}\label{thm:consistency}
Suppose that the observed data $\{Y_j:j=1,\ldots,n\}$ follow the model described
by equations (\ref{eq:basic}) and (\ref{eq:data_model}) and that assumptions
{\bf A1}--{\bf A4} are satisfied with $p\geq 3$. Suppose further that the
sequence $M$ is such that
\begin{equation}\label{eq:M_alpha_optimal}
c_1' \left(\frac{n}{\sigma_\varepsilon^2}\right)^{1/(2p+3)} \leq M \ll \left(\frac{n}{\sigma_\varepsilon^2\log n}\right)^{1/7}
\end{equation}
for some $c_1' > 0$,  $M^{-3/2}\ll \eta_M \ll M^{-1}$, and $\xi_n$ be
as defined in Lemma \ref{lem:consistency_x_j_hat}. Let $\bar\alpha_n := c_2' M^{-2}$ for some $c_2' > 0$
(sufficiently small) and
\begin{equation}\label{eq:alpha_n_def}
\alpha_n := C_0 M \max\left\{\sigma_\varepsilon \sqrt{\frac{M\log n}{n}}, M^{-(p+1)},\xi_n\right\},
\end{equation}
for some $C_0 > 0$.
Then as $n\to \infty$, with probability tending to one, there exists a local minimum $\widehat{\bs\beta}$
of the objective function $\tilde L_\delta(\bs\beta)$ (defined through (\ref{eq:loss_modified})),
which is also a global minimum within radius $\bar\alpha_n$ of $\bs\beta^*$ (note that, $\alpha_n \leq \bar\alpha_n$ by (\ref{eq:M_alpha_optimal})) such that, with $\widehat g := g_{\widehat{\bs\beta}}$,
\begin{equation}\label{eq:consistency}
\int_{x_{0,\delta}}^{x_{1,\delta}} |\widehat g(u) - g(u)|^2 du = O(\alpha_n^2).
\end{equation}
\end{theorem}
The proof of Theorem \ref{thm:consistency} is given in Section \ref{sec:proofs}.

\begin{remark}\label{rem:consistency}
Assuming $\sigma_\varepsilon$
to be a constant, if $M$ is chosen to be of the order $n^{1/(2p+3)}$, then
$\alpha_n^2$ in (\ref{eq:consistency}) simplifies to $n^{-2p/(2p+3)}\log n$, which is within
a factor of $\log n$ of the optimal rate in terms of the $L^2$-loss for estimating
$X'(t)$ based on the data $\{Y_j:j=1,\ldots,n\}$ given by (\ref{eq:basic}) when $X \in C^{p+1}([0,1])$.
The fact that an estimator of $g$ can attain this rate can be anticipated from the representation
of $g$ as $g = X' ~o~ X^{-1}$. For $p \geq 4$, we can improve the rate of convergence
of $\widehat g$ slightly further, by dropping the factor of $\log n$, as stated in the following result.
\end{remark}

\begin{theorem}\label{thm:optimal_rate}
Suppose that the conditions of Theorem \ref{thm:consistency}
are satisfied with $p \geq 4$ and, further, the sequence $M$ satisfies the condition
that $c (n/\sigma_\varepsilon^2)^{1/(2p+3)} \leq M \ll (n/\sigma_\varepsilon^2\log n)^{1/9}$
for some $c > 0$. Let $\widehat g$ be as in Theorem  \ref{thm:consistency}.
Then,
\begin{eqnarray}\label{eq:optimal_rate}
&& \int_{x_{0,\delta}}^{x_{1,\delta}} (\widehat g(x) - g(x))^2 dx \nonumber\\
&=& O_P\left(\frac{\sigma_\varepsilon^2 M^{3}}{n} \right)
+ O_P(M^{-2p}) + O_P\left(M^2 (\sigma_\varepsilon^2/n)^{2(p+1)/(2p+3)}\right), 
\end{eqnarray}
with the optimal rate given by $O_P((\sigma_\varepsilon^2/n)^{2p/(2p+3)})$, which is obtained when
$M = c (n/\sigma_\varepsilon^2)^{1/(2p+3)}$ for some $c > 0$.
\end{theorem}
Proof of Theorem \ref{thm:optimal_rate} is given in Section \ref{subsec:optimal_rate} of SM.

%\begin{remark}\label{rem:optimal_L_2_rate}
%We observe that the optimal rate of the $L^2$ norm of the error for estimating $g$ as given in
%Theorem \ref{thm:optimal_rate} is the optimal rate for nonparametric estimation of the derivative
%of an unknown smooth function with $p+1$ continuous derivatives.
%\end{remark}
%\begin{remark}\label{rem:rate_cubic_Bspline}
%Theorem \ref{thm:optimal_rate} shows that when $g \in C^4(D)$ (i.e.,
%$p=4$) and cubic B-splines with equally spaced knots are used for representing $g$, the best
%possible rate of convergence in $L^2$-loss is of the order $O_P(\alpha_n^2) \asymp O_P(n^{-8/11})$ and
%this is attained if $M \asymp n^{1/11}$ (again, assuming $\sigma_\varepsilon$ to be constant).
%\end{remark}

We can also  derive an approximate expression for the asymptotic variance of $\widehat g(\cdot)$.
Using a consistent root $\widehat{\bs\beta}$, we can use the equation $\nabla
L(\bs\beta)|_{\bs\beta = \widehat{\bs\beta}} = 0$.
%and apply a Taylor series
%expansion to write $\widehat{\bs\beta} = \bs\beta^* - [\nabla^2
%L(\bs\beta^*)]^{-1} \nabla L(\bs\beta^*) + O_P(\parallel \bs\beta -
%\bs\beta^*\parallel^2 + \parallel g - g_{\bs\beta^*}\parallel_{L^2}^2)$.
%Replacing the sample averages by corresponding population means conditional on $\mathbf{T}$,
Using the asymptotic representation of $\widehat{\bs\beta}-\bs\beta^*$ used
in the proof of Theorem \ref{thm:optimal_rate} (see Section \ref{subsec:optimal_rate} of SM),
and ignoring higher order terms
and the contribution of the model bias, and finally evaluating the expressions at
$\widehat{\bs\beta}$ instead of $\bs\beta^*$ (which is unknown), we have
\begin{equation}\label{eq:approx_asymptotic_variance}
\mbox{Var}(\widehat{\bs\beta}) \approx D(\widehat{\bs\beta}) :=
\widehat\sigma_\varepsilon^2 \left[ \sum_{j=1}^n \left(\frac{\partial
X(T_j;\widehat{\bs\beta})}{\partial \bs\beta}\right)\left(\frac{\partial
X(T_j;\widehat{\bs\beta})}{\partial \bs\beta}\right)^T\right]^{-1}.
\end{equation}
Here the estimated noise variance $\widehat\sigma_\varepsilon^2$ can be
computed as the mean squared error $(n-M)^{-1} \sum_{j=1}^M(Y_j -
X(T_j;\widehat{\bs\beta}))^2$. The expression
(\ref{eq:approx_asymptotic_variance}) allows us to obtain an approximate
asymptotic variance for $\widehat g(x)$ by $V(x) :=
\bs\phi(x)^T D(\widehat{\bs\beta}) \bs\phi(x)$, for any given $x$, where $\bs\phi(x) =
(\phi_{1,M}(x),\ldots,\phi_{M,M}(x))^T$.

\subsection{Initial estimator}\label{subsec:two_stage_regerssion}

In Theorem \ref{thm:consistency} we prove the rate of convergence
for a local minimizer, which is a global minimizer within a radius of $O(M^{-2})$ of $\bs\beta^*$
for a suitable range of values of $M$.
Therefore, we need an initial estimate which resides within this domain. In the following,
we describe one way of obtaining such an initial estimate, through a two-stage approach,
which is similar in spirit to the approaches by Chen and Wu (2008a, 2008b).
%$\widehat g$ derived in Theorem \ref{thm:consistency}.
%However, the rate of convergence of this estimator is faster than $O(M^{-2}) = O(\bar\alpha_n)$ if
%$M^4 \ll n^{(p+1)/(2p+3)}/\sqrt{\log n}$. So, for these range of $M$, which includes $M_*$,
%the two-stage estimator resides within the ball of radius $O(\bar\alpha_n)$ around $\bs\beta^*$,
%over which $\widehat g$, the optimizer of (\ref{eq:loss_modified}), is a global optimum. Therefore,
%treating $\widetilde g$ as initial value, with appropriate choice of $M$, largely ensures fast
%convergence to solution to $\widehat g$.

Suppose that we first estimate
$X(t)$ and $X'(t)$ by local polynomial smoothing  and denote these estimates by
$\hat X(t)$ and $\hat X'(t)$. Then, we fit the regression model
\begin{equation}\label{eq:two_stage_regression_model}
\hat X'(T_j) = \bs\phi(\hat X(T_j))^T \bs\beta + e_j, \qquad j=1,\ldots,n
\end{equation}
by ordinary least squares, where $\bs\phi = (\phi_{1,M},\ldots,\phi_{M,M})$.
We refer to the resulting estimator $\widetilde{\bs\beta}$ as the two-stage
estimator of $\bs\beta$:
\begin{equation}\label{eq:beta_tilde_two_stage}
\widetilde{\bs\beta} = [\sum_{j=1}^n \bs\phi(\hat X(T_j))\bs\phi(\hat X(T_j))^T \mathbf{1}_{[\delta,1-\delta]}(T_j)]^{-1}
(\sum_{j=1}^n \hat X'(T_j) \bs\phi(\hat X(T_j)) \mathbf{1}_{[\delta,1-\delta]}(T_j)).
\end{equation}
Since $X(t)$ is $p+1$ times continuously differentiable and $X'(t)$ is
$p$-times continuously differentiable (by {\bf A1}), and $\{\varepsilon_j\}$ is sub-Gaussian,
with the optimal choice of bandwidths, we have
\begin{eqnarray}
\max_{1\leq j \leq n} |\hat X(T_j) - X(T_j)| \mathbf{1}_{[\delta,1-\delta]}(T_j) &=&
O((\sigma_\varepsilon^2/n)^{(p+1)/(2p+3)}\sqrt{\log n} ) \label{eq:rate_local_polynomial}\\
\max_{1\leq j \leq n} |\hat X'(T_j) - X'(T_j)| \mathbf{1}_{[\delta,1-\delta]}(T_j)
&=& O((\sigma_\varepsilon^2/n)^{p/(2p+3)}\sqrt{\log n} ) \label{eq:rate_local_polynomial_deriv}
\end{eqnarray}
with probability tending to 1.
We state the following result about the rate of convergence of the two-stage estimator. The proof is given in
Section \ref{subsec:proofs_two_stage} of SM.
\begin{proposition}\label{prop:two_stage_regression}
Suppose that $p \geq 2$ and {\bf A1}--{\bf A4} hold and that the two-stage estimate of $g$ is given by
$\widetilde g := g_{\widetilde{\bs\beta}}$ where $\widetilde{\bs\beta}$ is defined in (\ref{eq:beta_tilde_two_stage}).
Then, supposing that $n^{1/(4p+6)} \ll M \ll n^{(p+1)/(4p+6)}/\sqrt{\log n}$, with probability tending to 1,
\begin{equation}\label{eq:two_stage_regression_rate}
\int_{x_{0,\delta}}^{x_{1,\delta}} |\widetilde g(u) - g(u)|^2 du = O(\widetilde \alpha_n^2)
\end{equation}
where
\begin{equation}
\widetilde \alpha_n = \max\{M^2(\sigma_\varepsilon^2/n)^{(p+1)/(2p+3)}\sqrt{\log n}, M^{-p}\}.
\end{equation}
\end{proposition}
When $\sigma_\varepsilon \asymp 1$, the optimal value of
$\widetilde \alpha_n$  is of the order $n^{-p(p+1)/(p+2)(2p+3)}(\log n)^{-p/(2p+4)}$
is obtained when $M \asymp M_* = n^{(p+1)/(p+2)(2p+3)}(\log n)^{-1/(2p+4)}$. It can be checked
that for all $p \geq 3$, this rate is slower than the optimal $\alpha_n$ for the
nonlinear regression-based estimator $\widehat g$ derived in Theorem \ref{thm:consistency}.
However, the rate of convergence of this estimator is faster than $O(M^{-2}) = O(\bar\alpha_n)$ if
$M^4 \ll n^{(p+1)/(2p+3)}/\sqrt{\log n}$. So, for these range of $M$, which includes $M_*$,
the two-stage estimator resides within the ball of radius $O(\bar\alpha_n)$ around $\bs\beta^*$,
over which $\widehat g$, the optimizer of (\ref{eq:loss_modified}), is a global optimum.

%Therefore, treating $\widetilde g$ as initial value, with appropriate choice of $M$, largely ensures fast
%convergence to solution to $\widehat g$.

\section{Proofs}\label{sec:proofs}

In this section, we outline the main steps of the proof. Some technical details are deferred to the Appendix.

The main idea behind the proof of Theorem \ref{thm:consistency} is
to obtain a lower bound on the difference $n^{-1} (L_\delta(\bs\beta) - L_\delta(\bs\beta^*))$ which
is proportional $\parallel \bs\beta - \bs\beta^*\parallel^2$ when $\bs\beta$ lies in an
annular region around $\bs\beta^*$. The outer
radius of the annular region depends on the degree of ill-conditioning of the problem,
as quantified by Proposition \ref{prop:kappa}, and the smoothness of the function $g$ and the approximating
bases, as indicated in condition {\bf A2}. This lower bound then naturally leads to the conclusion about
the existence and rate of convergence of a local minimizer $\widehat g$.

\subsection*{Proof of Proposition \ref{prop:kappa}}

For convenience of notations, we define $X_*(t)$ to be the sample path
$X(t;\bs{\beta}^*)$. Since $X^{\bs\beta}(\cdot;\bs\beta)$ is given by
(\ref{eq:tilde_X_beta_closed}) in the Appendix,
\begin{eqnarray*}
% X_{il}^{\theta_i}(t) &=&  e^{\theta_i} t
%g_{\bs{\beta}}( X_{il}(t)). \label{eq:tilde_X_theta_i_closed}\\
X^{\beta_r}(t) &=& g_{\bs{\beta}}( X(t))
\int_{ x_0}^{ X(t)} \frac{\phi_r(x)}{(g_{\bs{\beta}}(x))^2} dx, ~~r=1,\cdots, M,
%\label{eq:tilde_X_beta_closed}
\end{eqnarray*}
in order to prove Proposition
\ref{prop:kappa}, it suffices to find a lower bound on
\begin{equation*}
\min_{\parallel \mathbf{b}\parallel =1} \int_\delta^{1-\delta}
\left[\int_{\delta}^{t}g_{\mathbf{b}}(X_*(u))/g_{\bs{\beta}^*}(X_*(u))du
\right]^{2}\tilde f_T(t) dt
\end{equation*}
where $g_{\mathbf{b}}(u) = \mathbf{b}^T \bs{\phi}(u)$ with $\bs\phi =
(\phi_1,\ldots,\phi_M)^T$. By {\bf A3}, without loss of generality, we can take
the density $\tilde f_T(\cdot)$ to be uniform on $[\delta,1-\delta]$.

We make use of the following result known as Halperin-Pitt inequality
(Mitrinovic \textit{et al.}, 1991).
\begin{lemma}\label{lemma:Halperin_Pitt}
If $f$ is locally absolutely continuous and $f^{\prime\prime}$ is in
$L_{2}([0,A])$, then for any $\epsilon>0$ the following inequality holds with
$K(\epsilon)=1/\epsilon+12/A^{2}$,
\begin{equation}\label{eq:Halperin_Pitt}
\int_{0}^{A}(f^{\prime}(t))^2dt\leq K(\epsilon)\int_{0}^{A}f^{2}(t) dt
+\epsilon \int_{0}^{A}(f^{\prime\prime}(t))^2dt.
\end{equation}
\end{lemma}

Now defining
$$
R(t) :=
\int_{\delta}^{t}\frac{g_{\mathbf{b}}(X_*(u))}{g_{\bs{\beta}^*}(X_*(u))}du,
$$
%\end{equation*}
we have,
\begin{eqnarray*}
R'(t) &:=& \frac{dR(t)}{d t}
= \frac{g_{\mathbf{b}}(X_*(t))}{g_{\bs{\beta}^*}(X_*(t))} \\
R''(t) &:=& \frac{d^2R(t)}{dt^2} =
\left[\frac{g_{\mathbf{b}}'(X_*(t))}{g_{\bs{\beta}^*}(X_*(t))} -
\frac{g_{\mathbf{b}}(X_*(t))g_{\bs{\beta}^*}'(X_*(t))}{g_{\bs{\beta}^*}^2(X_*(t))}\right]
X_*'(t)\\
&=& \left[\frac{g_{\mathbf{b}}'(X_*(t))}{g_{\bs{\beta}^*}(X_*(t))} -
\frac{g_{\mathbf{b}}(X_*(t))g_{\bs{\beta}^*}'(X_*(t))}{g_{\bs{\beta}^*}^2(X_*(t))}\right]
g_{\bs{\beta}^*}(X_*(t)).
\end{eqnarray*}

%For sufficiently large $M$, by (\ref{eq:X_path_bias}) in the Appendix:
%\begin{equation*}
%\label{eq:X_path_bias}
%\parallel X(\cdot;\bs{\beta}^*) - X_g(\cdot)\parallel_\infty
%= O(M^{-p}),
%\end{equation*}

By (vi) of {\bf A2}, we have $\sup_{t \in [\delta,1-\delta]} |X_g(t) - X_*(t)| = O(M^{-(p+1)})$
and hence
\begin{equation}\label{eq:X_star_boundary}
X_*(1-\delta) \leq x_{1,\delta}  + |X_*(1-\delta)-x_{1,\delta}| < x_{1,M},
~ X_*(\delta) \geq x_{0,\delta} - |X_*(\delta)-x_{0,\delta}| >
x_{0,M}.
\end{equation}
Hence, using the facts that the coordinates of $\bs{\phi}(u)$ are $O(M^{1/2})$ and
the coordinates of $\bs{\phi}'(u)$ are $O(M^{3/2})$, and all these
functions are supported on intervals of length $O(M^{-1})$, we deduce that,
\begin{equation}\label{eq:r_prime_int}
\int_\delta^{1-\delta} (R''(t))^2  dt = O(M^2).
\end{equation}
An application of Lemma \ref{lemma:Halperin_Pitt} with $f(t) =  R(t-\delta)$
and $A=1-2\delta$ yields
\begin{equation}\label{eq:HP_main}
\int_\delta^{1-\delta} (R'(t))^2 dt \leq(1/\epsilon+12/(1-2\delta)^2)
\int_\delta^{1-\delta} (R(t))^2 dt  + \epsilon \int_\delta^{1-\delta}
(R''(t))^2 dt.
\end{equation}
Take $\epsilon = k_0 M^{-2}$ for some $k_0 > 0$, then by
(\ref{eq:r_prime_int}),
\begin{eqnarray*}
\int_{\delta}^{1-\delta} (R(t))^2 dt  &\geq& k_1 M^{-2}
\int_{\delta}^{1-\delta} (R'(t))^2 dt - k_2 M^{-2},
\end{eqnarray*}
for constants $k_1,k_2 > 0$ dependent on $k_0$. Next, we write
\begin{equation}\label{eq:smallest_eigen_bound}
\int_\delta^{1-\delta} (R'(t))^2 dt = \int_{X_*(\delta)}^{X_*(1-\delta)}
\frac{g_{\mathbf{b}}^2(v)}{g_{\bs{\beta}^*}^3(v)} dv =
\int_{X_*(\delta)}^{X_*(1-\delta)} g_{\mathbf{b}}^2(v) h(v) dv
\end{equation}
where $h(v) = g_{\bs{\beta}^*}^{-3}(v)$ which is bounded below by a positive
constant on the interval $[X_*(\delta),X_*(1-\delta)]$.
%\begin{equation}\label{eq:h_v}
%h(v) = g_{\bs{\beta}^*}^{-3}(v) \int \mathbf{1}_{\{x \leq v
%\leq X_1(1,x)\}} dF_a(x).
%\end{equation}
%If the knots are equally spaced on $[x_0
%+ \delta, x_1 - \delta]$ for some constant $\delta > 0$  bounded
%below, then $\inf_{v \in D_0} h(v)$ is bounded below (even as $M \to
%\infty$) where $D_0$ is the

Observe that by (\ref{eq:X_star_boundary}), the combined support of $\{\phi_{k,M}\}_{k=1}^M$, viz.,
$[x_{0,M},x_{1,M}]$,  contains (for sufficiently large $M$) the interval
$[X_*(\delta), X_*(1-\delta)]$. Also, $|x_{1,M} -
X_*(1-\delta)| \leq |x_{1,M} - x_{1,\delta}| + |x_{1,\delta} - X_*(1-\delta)| =
o(M^{-1})$ and $|x_{0,M} - X_*(\delta)| \leq |x_{0,M} - x_{0,\delta}| +
|x_{0,\delta} - X_*(\delta)| = o(M^{-1})$. These two facts and the condition
(v) of {\bf A2} imply that
\begin{eqnarray*}
&& \int_{X_*(\delta)}^{X_*(1-\delta)} g_{\mathbf{b}}^2(v)h(v) dv \\
&\geq&
\left(\inf_{v \in [X_*(\delta),X_*(1-\delta)]} h(v)\right)~ \mathbf{b}^T
[\int_{x_{0,M}}^{x_{1,M}} \bs\phi(v) (\bs\phi(v))^T dv - o(1)]\mathbf{b} ~\geq~ k_3,
\end{eqnarray*}
for some constant $k_3 > 0$, for sufficiently large $M$. Thus, by appropriate
choice of $\epsilon$, we have $\int_{\delta}^{1-\delta} (R(t))^2 dt \geq k_4
M^{-2}$ for some constant $k_4 > 0$, which yields (\ref{eq:G_star_condition}).

%If instead the knots are equally spaced on $[x_0,x_1]$ (excluding the boundaries),
%we have $\inf_{v \in D_0} h(v) \sim M^{-1}$ and hence, by choosing $\epsilon
%\sim M^{-3}$ we can show that $\int_0^1 \int_{0}^{t} (R(t,x))^2 dt dF_a(x) \geq k_5
%M^{-3}$, for some $k_5 > 0$, and consequently $\kappa_M = O(M^3)$.

\subsection*{Proof of Theorem \ref{thm:consistency}}

Define
\begin{equation}\label{eq:Gamma_beta_beta_star_def}
\Gamma_n(\bs\beta,\bs\beta^*) = \frac{1}{n}\sum_{j=1}^n (X(T_j;\bs\beta)
- X(T_j;\bs\beta^*))^2\mathbf{1}_{[\delta,1-\delta]}(T_j),
\end{equation}
${\cal A}_M(\alpha_n,\bar\alpha_n) = \{\bs\beta \in \mathbb{R}^M : \alpha_n \leq \parallel \bs\beta - \bs\beta^*\parallel \leq \bar\alpha_n\}$,
and
$$
D_n^*=  \frac{1}{n}\sum_{j=1}^n (X_g(T_j) - X(T_j;\bs\beta^*))^2\mathbf{1}_{[\delta,1-\delta]}(T_j).
$$
Suppose that $\bs\beta \in {\cal A}_M(\alpha_n,\bar\alpha_n)$. Henceforth, we use $X(t;\bs\beta)$ to denote $X(t;\bs\beta;x_{0,\delta})$
and $X_g(t)$ to denote $X_g(t;x_{0,\delta})$. Then
\begin{eqnarray}\label{eq:diff_L_delta_beta}
&&~~\frac{1}{n} L_\delta(\bs\beta) - \frac{1}{n} L_\delta(\bs\beta^*) \nonumber\\
&=& \frac{1}{n} \sum_{j=1}^n (Y_j - X(T_j;\bs\beta))^2
\mathbf{1}_{[\delta,1-\delta]}(T_j) - \frac{1}{n} \sum_{j=1}^n (Y_j - X(T_j;\bs\beta^*))^2
\mathbf{1}_{[\delta,1-\delta]}(T_j) \nonumber\\
&=& \frac{1}{n}\sum_{j=1}^n (X(T_j;\bs\beta) - X(T_j;\bs\beta^*))^2\mathbf{1}_{[\delta,1-\delta]}(T_j) \nonumber\\
&& ~~~- \frac{2}{n} \sum_{j=1}^n \varepsilon_j (X(T_j;\bs\beta) - X(T_j;\bs\beta^*))\mathbf{1}_{[\delta,1-\delta]}(T_j) \nonumber\\
&& ~~~- \frac{2}{n} \sum_{j=1}^n (X_g(T_j) - X(T_j;\bs\beta^*))(X(T_j;\bs\beta) - X(T_j;\bs\beta^*))\mathbf{1}_{[\delta,1-\delta]}(T_j),
\end{eqnarray}
where $U_{1n}(\bs\beta,\bs\beta^*)$ and  $U_{2n}(\bs\beta,\bs\beta^*)$, are the second and third summations in
the above expression, respectively.
Next, we write
\begin{eqnarray}\label{eq:diff_L_L_tilde_delta_beta}
&&~~\frac{1}{n} \tilde L_\delta(\bs\beta) - \frac{1}{n}  L_\delta(\bs\beta) \nonumber\\
&=& \frac{1}{n} \sum_{j=1}^n \left[(Y_j -  X(T_j;\bs\beta,\widehat x_0))^2
- (Y_j -  X(T_j;\bs\beta))^2\right]\mathbf{1}_{[\delta,1-\delta]}(T_j)
\nonumber\\
&=& \frac{1}{n}\sum_{j=1}^n (X(T_j;\bs\beta;\widehat x_0)-X(T_j;\bs\beta))^2  \mathbf{1}_{[\delta,1-\delta]}(T_j) \nonumber\\
&& ~~~-\frac{2}{n}\sum_{j=1}^n \varepsilon_j (X(T_j;\bs\beta;\widehat x_0)-X(T_j;\bs\beta)) \mathbf{1}_{[\delta,1-\delta]}(T_j)
\nonumber\\
&& ~~~- \frac{2}{n} \sum_{j=1}^n (X(T_j;\bs\beta;\widehat x_0)-X(T_j;\bs\beta))(X_g(T_j) - X(T_j;\bs\beta)) \mathbf{1}_{[\delta,1-\delta]}(T_j),
\end{eqnarray}
where $V_{1n}(\bs\beta)$,  $V_{2n}(\bs\beta)$, $V_{3n}(\bs\beta)$ are the three summations in the last
expression.

From (\ref{eq:diff_L_delta_beta}) and (\ref{eq:diff_L_L_tilde_delta_beta}), we deduce that
\begin{eqnarray}\label{eq:diff_L_tilde_delta_beta_expand}
&& ~~\frac{1}{n} \tilde L_\delta(\bs\beta) - \frac{1}{n} \tilde L_\delta(\bs\beta^*) \nonumber\\
&=&
\frac{1}{n} \left(L_\delta(\bs\beta) - L_\delta(\bs\beta^*)\right)
+ \frac{1}{n} \left(\tilde L_\delta(\bs\beta) - L_\delta(\bs\beta)\right)
- \frac{1}{n} \left(\tilde L_\delta(\bs\beta^*) - L_\delta(\bs\beta^*)\right)  \nonumber\\
&=& \Gamma_n(\bs\beta,\bs\beta^*) -U_{1n}(\bs\beta,\bs\beta^*) - U_{2n}(\bs\beta,\bs\beta^*)
+ (V_{1n}(\bs\beta)-V_{1n}(\bs\beta^*))  \nonumber\\
&&~ - (V_{2n}(\bs\beta)-V_{2n}(\bs\beta^*))  + U_{3n}(\bs\beta,\bs\beta^*) - U_{4n}(\bs\beta,\bs\beta^*) 
\end{eqnarray}
where
\begin{eqnarray*}
U_{3n}(\bs\beta,\bs\beta^*) &\hskip-.1in=& \hskip-.1in\frac{2}{n} \sum_{j=1}^n (X_g(T_j) - X(T_j;\bs\beta^*))
(X(T_j;\bs\beta;\widehat x_0) - X(T_j;\bs\beta)))\mathbf{1}_{[\delta,1-\delta]}(T_j) \\
&& \hskip-.2in -\frac{2}{n} \sum_{j=1}^n (X_g(T_j) - X(T_j;\bs\beta^*))(X(T_j;\bs\beta^*;\widehat x_0) - X(T_j;\bs\beta^*))
\mathbf{1}_{[\delta,1-\delta]}(T_j)  \nonumber\\
U_{4n}(\bs\beta,\bs\beta^*) &\hskip-.1in=& \hskip-.1in\frac{2}{n} \sum_{j=1}^n (X(T_j;\bs\beta)
- X(T_j;\bs\beta^*))(X(T_j;\bs\beta;\widehat x_0)-X(T_j;\bs\beta))
\mathbf{1}_{[\delta,1-\delta]}(T_j).
\end{eqnarray*}
Using the fact that $\equiv X^a(t;\bs\beta,a_0) := (\partial/\partial a) X(t;\bs\beta,a)|_{a=a_0}$
can be expressed as
\begin{equation*}
X^a(t;\bs\beta,a_0) = \frac{g_{\bs\beta}(X(t;\bs\beta,a_0))}{g_{\bs\beta}(a_0)}, ~~~t\in [\delta,1-\delta],
\end{equation*}
provided $g_{\bs\beta}(x) > 0$ for $x \in [a_0, X(1-\delta;\bs\beta,a_0)]$, we have, for all
$\bs\beta \in {\cal A}_M(\alpha_n,\bar\alpha_n)$,
\begin{equation}\label{eq:X_diff_partial_a_bound}
\sup_{a_0 \in [x_{0,\delta} - \xi_n,x_{0,\delta}+\xi_n]} \sup_{t \in [\delta,1-\delta]} |X(t;\bs\beta,a_0) - X(t;\bs\beta,x_{0,\delta})|
\leq C_1\xi_n
\end{equation}
for some $C_1 > 0$. Here, we have used the fact that
for $t \in [\delta,1-\delta]$, and $a_0 \in [x_{0,\delta} - \xi_n,x_{0,\delta}+\xi_n]$,
\begin{equation}
X(t;\bs\beta,a_0) = \tilde G_{\bs\beta}^{-1}(t-\delta + G_{\bs\beta}(a_0))
~~\mbox{where}~~ \tilde G_{\bs\beta}(x) := \int_{x_{0,M}}^x \frac{du}{g_{\bs\beta}(u)}~,
\end{equation}
and that
$$
\sup_{\bs\beta \in {\cal A}(\alpha_n,\bar\alpha_n)}
\sup_{x\in [x_{0,M},x_{1,M}]}| g_{\bs\beta}(x) - g_{\bs\beta^*}(x)| = O(\bar\alpha_n M^{1/2}) = O(M^{-3/2}),
$$
so that, by using (\ref{eq:X_diff_bound}), and the fact that $M^{-3/2} \ll \eta_M \ll M^{-1}$,
$$
[x_{0,\delta} - \xi_n,\sup_{\bs\beta \in {\cal A}(\alpha_n,\bar\alpha_n)}\sup_{a_0 \in [x_{0,\delta} - \xi_n,x_{0,\delta}+\xi_n]} X(1-\delta;\bs\beta,a_0)]  \subset [x_{0,M},x_{1,M}]
$$
for large enough $M$ and $n$.

We now bound individual terms in the expansion (\ref{eq:diff_L_tilde_delta_beta_expand}). First,
we have the following lower bound on $\Gamma_n(\bs\beta,\bs\beta^*)$, the proof of which is given in Appendix C.
\begin{lemma}\label{lemma:Gamma_beta_beta_star_bound}
Let $\Gamma_n(\bs\beta,\bs\beta^*)$ be as defined in (\ref{eq:Gamma_beta_beta_star_def}). Then given $\eta > 0$, there exist
constants $d_1(\eta) > 0$ and $d_2,d_3 > 0$ independent of $\eta$ such that
\begin{equation}\label{eq:Gamma_beta_beta_star_bound}
\Gamma_n(\bs\beta,\bs\beta^*) \geq d_1(\eta) \frac{1}{M^2} \parallel \bs\beta - \bs\beta^*\parallel^2
-  d_2 \parallel \bs\beta - \bs\beta^*\parallel^4 M^2(1+ d_3\bar\alpha_n^2 M^3)
\end{equation}
uniformly in $\bs\beta \in {\cal A}_M(\alpha_n,\bar\alpha_n)$ with probability at least $1-n^{-\eta}$.
\end{lemma}
Since $\bar\alpha_n M^{3/2} = c_2'M^{-1/2} = o(1)$, and the constant $c_2'$
can be chosen to be small enough so that we can conclude from (\ref{eq:Gamma_beta_beta_star_bound}) that given
$\eta > 0$, there exists $d_4(\eta) > 0$ such that
\begin{equation}\label{eq:Gamma_beta_beta_star_bound_refined}
\mathbb{P}\left(\Gamma_n(\bs\beta,\bs\beta^*) \geq \frac{d_4(\eta)}{M^2} \parallel \bs\beta - \bs\beta^*\parallel^2
~~\mbox{for all}~\bs\beta \in {\cal A}_M(\alpha_n,\bar\alpha_n)\right) \geq 1- n^{-\eta}.
\end{equation}
Next, by Cauchy-Schwarz inequality, we have
\begin{equation}\label{eq:U_2n_bound}
|U_{2n}(\bs\beta,\bs\beta^*)| \leq 2\sqrt{D_n^*} \sqrt{\Gamma_n(\bs\beta,\bs\beta^*)}.
\end{equation}
Next, by (\ref{eq:X_diff_partial_a_bound}), we have
\begin{equation}\label{eq:V_1n_bound}
\max\{V_{1n}(\bs\beta^*),\sup_{\bs\beta \in {\cal A}(\alpha_n,\bar\alpha_n)}V_{1n}(\bs\beta)\}
\leq C_1^2 \xi_n^2,
\end{equation}
and hence
\begin{equation}\label{eq:U_3n_bound}
\sup_{\bs\beta \in {\cal A}(\alpha_n,\bar\alpha_n)}|U_{3n}(\bs\beta,\bs\beta^*)| \leq 4C_1 \xi_n \sqrt{D_n^*}
\end{equation}
and
\begin{equation}\label{eq:U_4n_bound}
|U_{4n}(\bs\beta,\bs\beta^*)| \leq 2C_1\xi_n \sqrt{\Gamma_n(\bs\beta,\bs\beta^*)}.
\end{equation}
Next, defining
$$
Z(\bs\beta) = \frac{\sum_{j=1}^n \varepsilon_j (X(T_j;\bs\beta) - X(T_j;\bs\beta^*))\mathbf{1}_{[\delta,1-\delta]}(T_j)}
{\sigma_\varepsilon \sqrt{\sum_{j=1}^n (X(T_j;\bs\beta) - X(T_j;\bs\beta^*))^2\mathbf{1}_{[\delta,1-\delta]}(T_j)}}~,
$$
and setting $Z(\bs\beta)$ being zero if the denominator is zero, we have
\begin{equation}\label{eq:U_1n_bound}
|U_{1n}(\bs\beta)| \leq \frac{2\sigma_\varepsilon}{\sqrt{n}} \sqrt{\Gamma_n(\bs\beta,\bs\beta^*)} |Z(\bs\beta)|.
\end{equation}
Let ${\cal B}_M(\Delta;\alpha_n,\bar\alpha_n)$ be a $\Delta$-net for ${\cal A}_M(\alpha_n,\bar\alpha_n)$. Then
$|{\cal B}_M(\Delta;\alpha_n,\bar\alpha_n)| \leq 3(\bar\alpha_n/\Delta)^M$. Then, by using Lemma \ref{lem:subgaussian_Hoeffding} in SM,
and (\ref{eq:Gamma_beta_beta_star_bound_refined}),
we conclude that given $\eta > 0$, there exist constants $c_1(\eta) > 0, C'(\eta) > 0$, and a set $A_{1\eta}$ with
$\mathbb{P}(\mathbf{T} \in A_{1\eta}) \geq 1-n^{-\eta}$,
such that  for all $\mathbf{T} \in A_{1\eta}$,
\begin{eqnarray*}
&& \mathbb{P}\left(\max_{\bs\beta \in {\cal B}_M(\delta;\alpha_n,\bar\alpha_n)}|Z(\bs\beta)|
> c_1(\eta)\sqrt{M \log(\bar\alpha_n/\delta)}~|~ \mathbf{T}\right)
\leq C'(\eta) \left(\frac{\Delta}{\bar\alpha_n}\right)^{\eta M}
\end{eqnarray*}
for some constant $C' > 0$. Thus, taking $\delta$ to be sufficiently small, say,
$\delta = n^{-c}$ for $c$ large enough, and using the
smoothness of the process $Z(\bs\beta)$ as a function of $\bs\beta$,  we can show that
given any $\eta > 0$, there exists $c_2(\eta) > 0$, such that for all $\mathbf{T} \in A_{1\eta}$,
\begin{equation}\label{eq:Z_beta_sup_bound}
\mathbb{P}\left(\sup_{\bs\beta \in {\cal A}_M(\alpha_n,\bar\alpha_n)}
|Z(\bs\beta)|
\leq ~c_2(\eta)\sqrt{M \log n}~\left|\right.~\mathbf{T}\right) > 1-n^{-\eta}.
\end{equation}
Very similarly, defining
$$
\tilde  Z(\bs\beta) = \frac{\sum_{j=1}^n \varepsilon_j (X(T_j;\bs\beta;\widehat x_0)-X(T_j;\bs\beta))\mathbf{1}_{[\delta,1-\delta]}(T_j)}
{\sigma_\varepsilon \sqrt{\sum_{j=1}^n (X(T_j;\bs\beta;\widehat x_0)-X(T_j;\bs\beta))^2\mathbf{1}_{[\delta,1-\delta]}(T_j)}},
$$
expressing
$V_{2n}(\bs\beta) = 2\sigma_\varepsilon n^{-1/2} \sqrt{V_{1n}(\bs\beta)} \tilde Z(\bs\beta)$,
and using (\ref{eq:V_1n_bound}),
we have, for any given $\eta > 0$, there exists $c_3(\eta) > 0$ and a set $A_{2\eta}$ with
$\mathbb{P}(\mathbf{T} \in A_{2\eta}) \geq 1-n^{-\eta}$,
such that  for all $\mathbf{T} \in A_{2\eta}$,
\begin{equation}\label{eq:V_2n_bound}
\mathbb{P}\left(\sup_{\bs\beta \in {\cal A}_M(\alpha_n,\bar\alpha_n)\cup \{\bs\beta^*\}}
|V_{2n}(\bs\beta)|
\leq ~c_3(\eta)\sigma_\varepsilon \xi_n\sqrt{\frac{M \log n}{n}}~\left|\right.~\mathbf{T}\right) > 1-n^{-\eta}.
\end{equation}
Finally, by {\bf A2} we have the bound
\begin{equation}\label{eq:D_n_star_bound}
D_n^* \leq \sup_{t \in [\delta,1-\delta]} |X_g(t) - X(t;\bs\beta^*)|^2 \leq C_2 M^{-2(p+1)}
\end{equation}
for some $C_2 > 0$.

Combining (\ref{eq:U_2n_bound})--(\ref{eq:D_n_star_bound}),
we claim that, given $\eta > 0$, there exist constants $C_3(\eta) > 0$, $C_4(\eta) > 0$,
and constants $C_l > 0$, $l=5,\ldots,8$, not depending on $\eta$, such that
uniformly on ${\cal A}_M(\alpha_n,\bar\alpha_n)$
\begin{eqnarray}\label{eq:L_tilde_beta_diff_bound}
&& \frac{1}{n}\tilde L_\delta(\bs\beta) - \frac{1}{n}\tilde L_\delta(\bs\beta^*) \nonumber\\
&\geq& \Gamma_n(\bs\beta,\bs\beta^*)
- \sqrt{\Gamma_n(\bs\beta,\bs\beta^*)}\left(C_3(\eta) \sqrt{\frac{M\log n}{n}} + C_5 M^{-(p+1)} + C_6 \xi_n \right) \nonumber\\
&& ~~~
- \xi_n \left(C_4(\eta) \sqrt{\frac{M \log n}{n}} + C_7 M^{-(p+1)} + C_8\xi_n \right) 
\end{eqnarray}
with probability at least $1-O(n^{-\eta})$.

From (\ref{eq:L_tilde_beta_diff_bound}) and (\ref{eq:Gamma_beta_beta_star_bound_refined}), and a careful choice of the
constant $C_0$ in the definition (\ref{eq:alpha_n_def}) of $\alpha_n$, and with  $M$ as in \ref{eq:M_alpha_optimal},
we conclude that for any $\eta > 0$,
there exists $C_{9}(\eta) > 0$ such that, uniformly in $\bs\beta \in {\cal A}_M(\alpha_n,\bar\alpha_n)$,
\begin{equation}
\frac{1}{n}\tilde L_\delta(\bs\beta) - \frac{1}{n}\tilde L_\delta(\bs\beta^*) \geq C_{9}(\eta) \frac{1}{M^2} \parallel \bs\beta - \bs\beta^*\parallel^2
\end{equation}
with probability at least $1-O(n^{-\eta})$. From this, we can conclude that with probability at least $1-O(n^{-\eta})$
there exists a local minimum $\widehat{\bs\beta}$ of $\tilde L_\delta(\bs\beta)$, which is also a global minimum within radius $\bar\alpha_n$
of $\bs\beta^*$ and which satisfies $\parallel \widehat{\bs\beta} - \bs\beta^*\parallel = O(\alpha_n)$.

%These results depend  on the application of perturbation bounds for differential
%equations and appropriate deviation inequalities. Also, in the proofs of
%these lemmas (see the Appendix) and Theorem \ref{thm:consistency}, we shall
%repeatedly use the following key decomposition:
%\begin{equation}\label{eq:error_decomp}
%Y_j - X_j(\bs{\beta}) = \varepsilon_j + (X_j^g - X_j(\bs{\beta}^*)) +
%(X_j(\bs{\beta}^*) - X_j(\bs{\beta})).
%\end{equation}
%Under the assumptions of Theorem \ref{thm:consistency}, $\alpha_n M^{3/2}
%\to 0$ and thus the perturbation bounds (\ref{eq:X_path_bias}) -
%(\ref{eq:X_diff_beta_beta_bound}) are valid.

\section{Simulation Study}
\label{sec:simulation}

In this section, we conduct a simulation study to examine the finite sample
performance of the proposed estimation procedure, as well as to compare it with
the two-stage estimator described in Section \ref{subsec:two_stage_regerssion}.

In the simulation, the true gradient function $g$ is represented by $4$
B-spline functions with knots at $0.35, 0.60, 0.85, 1.10 $ and respective
coefficients $0.1, 1.2, 1.6, 0.4$ (shown by the blue curve in Figure
\ref{fig:simu}). We set the initial value $X(0)=x_0=0.25$ in equation
(\ref{eq:basic}) to generate the true trajectory $X(\cdot)$. We then simulate
$100$ independent data sets according to equation (\ref{eq:data_model}) .
Specifically,  for each data set, we first randomly choose an integer $n$ from
$\{60,\cdots, 100\}$. Then $n$ observation times $\{t_1,\cdots, t_n\}$ are
uniformly sampled from $[0,1]$. Finally, the $Y_j$'s are generated according to
equation (\ref{eq:data_model}) with added noise $\epsilon_i \sim {\rm
Normal}(0,0.01^2)$. The observed data from one such replicate is shown in
Figure \ref{fig:obs_trajectory} in SM together with the true trajectory $X(\cdot)$.

We fit the proposed estimator $\widehat{g}(\cdot)$ with $M$ B-spline basis functions with equally
spaced knots on $[0.1,1.1]$. We consider $M=3,4,5$ and choose $M$ by an
approximate leave-one-out CV score criterion similar to that used in Paul et al. (2011).
%\textcolor{blue}{More details on this approach can be found in Paul et al. (2011).}
Out of the $100$ replicates, $43$ times the model with $M=4$ (the true model) is chosen and $66$
times the model with $M=5$ is chosen.
%%% initial value for beta: 1,...,1. true x0 is used.

We also consider the two-stage estimator, where in the first stage, the sample trajectory $X(\cdot)$ and its
derivative $X'(\cdot)$ are estimated by applying local linear and local
quadratic smoothing with Gaussian Kernel, respectively, to the observed data
$\{(t_j, Y_j)\}_{j=1}^n $. The bandwidths are chosen by cross-validation. In
the second stage, a quadratic smoothing of  $\widehat{X}'(\cdot)$ versus
$\widehat{X}(\cdot)$  is performed to get an estimate of $g(\cdot)$.

Figure \ref{fig:simu} shows the estimated gradient functions (red curves) of
these $100$ independent replicates overlayed on the true gradient function
(blue curve). It can be seen from this figure that, the proposed estimator
shows little bias. Its sampling variability is somewhat larger  on the left side of
the observed $x$ domain than on the right side of the observed $x$ domain. It
performs much better than the two-stage estimator which shows both high bias
and high variance. Indeed, the bias of the two-stage estimator would not go
away even when in the second stage the true model is used to estimate $g$
(through a least-squares regression of $\widehat{X}'(\cdot)$ versus
$\widehat{X}(\cdot)$).

Figure \ref{fig:simu_traj} shows the estimated trajectories (red curves) of
these $100$ independent replicates overlayed on the true trajectory (blue
curve). In the left panel of the figure, the estimated trajectories are solved
from equation (\ref{eq:basic}) using the 4th-order Runge-Kutta method with $g$
being the proposed estimator $\widehat{g}(\cdot)$. In the right panel of the
figure, the trajectories are estimated by applying local linear smoothing of
the observed data (which are then used in the two-stage fitting for
$g(\cdot)$). The estimated trajectories from the proposed procedure follow the
true trajectory very well with little bias, whereas the estimator from the
first-stage smoothing of the two-stage procedure shows more bias and more
variability. Figure \ref{fig:simu_traj_deriv} in SM shows the estimated derivative of
the trajectory. Again, the proposed procedure gives a much better estimate of
$X^{\prime}(\cdot)$ than the presmoothing estimate (by local quadratic
smoothing) used in the two-stage procedure.

%%% First stage, Gauss kernel, CV to choose bandwith.

\section{Application : Berkeley Growth Data}
\label{sec:real}

In this section, we apply the proposed model to the Berkeley
growth data (Tuddenham and Snyder, 1954). Although in the literature, there are
many studies of growth curves (Hauspie et al., 1980; Milani, 2000), most of
them try to model either the growth trajectories (i.e., $X(\cdot)$) or the rate
of growth (i.e., $X'(\cdot)$). On the contrary, our goal is to estimate the
gradient function, i.e., the functional relationship between  $X'(\cdot)$ and
$X(\cdot)$ which provides insights of the growth dynamics, such as  at what height the growth rate tends to
be the highest.

Specifically, we fit the proposed model to  each of the $54$ female subjects in
this data set. For each girl, her heights were measured at $31$ time points
from $1$ year old to $18$ years old. We use $M$ B-spline basis functions with
equally spaced knots.
%starting from the shortest height in the data set ($\sim
%46$ cm) and ending with the tallest height in the data set ($\sim 183$ cm).
We consider  $M=4,5,6,7$ and for each subject we choose the ``best'' $M$ using an
approximate leave-one-out CV score. In $37$ out of $54$ subjects, the model
with $M=6$ is chosen, and for the rest $17$ subjects, the model with $M=7$ is
chosen. Figure \ref{fig:fit_gradient} shows the fitted gradient functions for
these $54$ subjects. From this figure, we can see that, most girls experienced
two growth spurs, one at the birth (when their heights are shortest) and
another when they were around either $130$ cm tall or $150$ cm tall. Moreover
Figure \ref{fig:fit_gradient_SE} in SM shows the fitted gradient functions with the
two-standard-error bands (by equation (\ref{eq:approx_asymptotic_variance}))
for $25$ girls. Figure \ref{fig:fit_trajectory} in SM shows the observed (red dots)
and fitted (black curve) growth trajectories for these $25$ girls. It can be
seen that, the fitted trajectories fit the observed data very well.

\section{Discussion}\label{sec:discussion}

In this paper we have proposed an estimation procedure for nonparametrically
estimating the unknown gradient function of a first order autonomous differential
equation over a finite domain, when the trajectories are strictly monotone. In this
section, we discuss the asymptotic rate optimality of the proposed estimator. We show
that, if the estimators of the gradient function $g$ are restricted to a class of uniformly
Lipschitz function, the optimal rate for estimation of $g$, i.e., of the order
$n^{-2p/(2p+3)}$, is the same as the optimal rate for
estimation of the derivative of $X \equiv X_g$ based on model (\ref{eq:data_model})
in terms of the $L^2$ loss.
We conjecture that the Lipschitz requirement on the estimator of $g$ is not necessary and
the minimax rate for estimation of $g$ is indeed of the order $n^{-2p/(2p+3)}$.

In order to make this statement precise, we first specify the function class for $g$ as
\begin{equation}\label{eq:g_class}
{\cal G} = \{g : D \to \mathbb{R}_+ : c_0 \leq g \leq c_1; |g'| \leq c_2; g \in C^p(D)\}
\end{equation}
where $0 < c_0 < c_1 <\infty$ and $0 < c_2 < \infty$ are constants. Define the class of
uniformly Lipschitz functions
\begin{equation*}
{\cal L} = \{h : D \to \mathbb{R} : |h(x) - h(y)| \leq c_4 |x-y|~\mbox{for all}~ x, y \in D\}
\end{equation*}
where $c_4 \in (0,\infty)$ depends on (at least as large as) $c_2$ in (\ref{eq:g_class}).
If $g \in {\cal G}$, then we have $X_g \in C^{p+1}([0,1])$ and $X_g' \in C^p([0,1])$.
In addition, we assume the observation model (\ref{eq:data_model}) with the noise
$\varepsilon_i \stackrel{i.i.d.}{\sim} N(0,\sigma_\varepsilon^2)$.

Let $\delta$ be as in Section \ref{sec:model}. By the condition $c_0 \leq g \leq c_1$,
we know that there exist $0 < c_0(\delta) < c_1(\delta) < \infty$ such that
$c_0(\delta) \leq X_g(t) \leq c_1(\delta)$
for all $t \in [\delta,1-\delta]$, for all $g \in {\cal G}$. Define,
$\parallel f \parallel_{2,\delta} = (\int_{\delta}^{1-\delta} (f(t))^2 dt)^{1/2}$.
Then there are constants $c_2(\delta), c_3(\delta) > 0$ such that
for any given estimator $\widehat g \in {\cal L}$ of $g$,
\begin{eqnarray}\label{eq:g_hat_error_bound}
c_2(\delta)\parallel \widehat g~o~X_g - g~o~X_g \parallel_{2,\delta}^2
&\leq& \int_{X_g(\delta)}^{X_g(1-\delta)} |\widehat g(u) - g(u)|^2 du \nonumber\\
&\leq& c_3(\delta) \parallel \widehat g~o~X_g - g~o~X_g \parallel_{2,\delta}^2. 
\end{eqnarray}
Observe that $g~o~X_g = X_g'$.

On the other hand, since $X_g \in C^{p+1}([0,1])$, there exists
an estimator $\widehat X_{op}$ with the property that,
given $\epsilon > 0$, there exists constant $K_1(\epsilon) > 0$ such that
\begin{equation}\label{eq:X_g_minimax}
\sup_{g \in {\cal G}}\mathbb{P}(\parallel \widehat X_{op} -
X_g \parallel_{2,\delta}^2 > K_1(\epsilon) n^{2(p+1)/(2p+3)}) < \epsilon
\end{equation}
for all $n \geq N_1(\epsilon)$.

We define the estimator $\widetilde X' := \widehat g~o~\widehat X_{op}$
for $X_g'$.  Then, by triangle inequality,
\begin{eqnarray}\label{eq:X_deiv_estimate_bound}
\parallel \widehat g~o~X_g - g~o~X_g \parallel_{2,\delta} &\hskip-.1in=&
\hskip-.1in\parallel \widehat g~o~X_g - X_g' \parallel_{2,\delta} \nonumber\\
&\hskip-.1in\geq&
\hskip-.1in\parallel \widetilde X' - X_g'\parallel_{2,\delta} -
\parallel \widehat g~o~\widehat X_{op} - \widehat g~o~X_g \parallel_{2,\delta} \nonumber\\
&\hskip-.1in\geq& \hskip-.1in\parallel \widetilde X' - X_g'\parallel_{2,\delta} - c_4 \parallel \widehat X_{op} -
X_g \parallel_{2,\delta}, 
\end{eqnarray}
where, in the last step we have used the fact that $\widehat g \in {\cal L}$.

Since $X_g' \in C^p([0,1])$, the minimax rate of estimation of $X_g'$
in terms of the $L^2$ loss $\parallel \cdot \parallel_{2,\delta}^2$ is of the
order $n^{-2p/(2p+3)}$. This can be derived directly for $g$ restricted to
${\cal G}$ by only slightly modifying the arguments in Stone (1982).
Combining this fact with (\ref{eq:g_hat_error_bound}),
(\ref{eq:X_g_minimax}) and (\ref{eq:X_deiv_estimate_bound}), we obtain that
there exists $K_2 > 0$, such that
%for any $\epsilon > 0$, there exists $K_2(\epsilon) > 0$ such that
\begin{equation*}
\lim\inf_{n\to \infty} \inf_{\widehat g \in {\cal L}} \sup_{g \in {\cal G}}
\mathbb{P}\left( \int_{X_g(\delta)}^{X_g(1-\delta)} |\widehat g(u) - g(u)|^2 du >  K_2 n^{-2p/(2p+3)}\right)
> 0.
\end{equation*}
%for all $n \geq N_2(\epsilon)$.
In other words, as long as $\widehat g$ is uniformly Lipschitz, the rate
$n^{-2p/(2p+3)}$ is a lower bound on the rate for estimating $g$ in terms of the $L^2$-loss. We note that, the requirement
$\widehat g \in {\cal L}$ can be relaxed by only requiring that this holds with probability approaching one
as $n \to \infty$. The latter is satisfied by the estimator we proposed.  Thus, combining with
Theorem \ref{thm:optimal_rate}, we deduce that the optimal rate of estimation of $g$ is $n^{-2p/(2p+3)}$ for
$p \geq 4$.

%%initial x0: first observation; beta_ini: from the projected empirical fitted gradient

%\clearpage

%%%%%%%%%%%%%%%%%%%%%%%%%%%%%%%%%%%%%%%%%%%%%%%%%%%%%%%
%%%%%%%%%% Appendix %%%%%%%%%%%%%%%%%%%%%
%%%%%%%%%%%%%%%%%%%%%%%%%%%%%%%%%%%%%%%%%%%%%%%%%%%%%%%

%\begin{center}
%\Large{\bf Supplementary Materials for ``Consistency in a semiparametric
%autonomous dynamical system'' }
%\end{center}

%\setcounter{page}{1}

\setcounter{equation}{0}
\renewcommand{\theequation}{A.\arabic{equation}}
\setcounter{figure}{0}
\renewcommand{\thefigure}{A.\arabic{figure}}
\setcounter{proposition}{0}
\renewcommand{\theproposition}{A.\arabic{proposition}}
\setcounter{lemma}{0}
\renewcommand{\thelemma}{A.\arabic{lemma}}
\setcounter{corollary}{0}
\renewcommand{\thecorollary}{A.\arabic{corollary}}

\section{Appendix}\label{sec:Appendix}

In this section, we provide technical details for the
proofs of the main results. Specifically, in Appendix A, we present results on
perturbation analysis of differential equations that are central to controlling
the bias in the estimates. In Appendix B, we verify that condition (vi) of {\bf A2}
is satisfied by a B-spline basis of sufficiently high order.
%In Appendix C, we list the relevant results on sub-Gaussian random variables.
In Appendix C, we prove Lemma \ref{lemma:Gamma_beta_beta_star_bound}. Further technical details
are given in the Supplementary Material.
%In Appendix E, we prove Theorem \ref{thm:optimal_rate}, and in Appendix F, we analyze the rate
%of convergence for the two stage estimator defined in Section \ref{subsec:two_stage_regerssion}.

\subsection*{Appendix A : Properties of sample trajectories and their
derivatives}\label{sec:Appendix_A}

%\subsection*{Derivatives of the sample paths $\{X_{il}(\cdot)\}$
%with respect to $(\boldsymbol{\theta},\boldsymbol{\beta})$}

Throughout this subsection, with slight abuse of notation, we use $X(\cdot)$ to mean $X(\cdot;\bs\beta)$,
unless otherwise noted.

Since $X(\cdot)$ satisfies the ODE
\begin{equation}\label{eq:tilde_X_i}
X(t) = x_0 + \int_{0}^t \sum_{k=1}^M \beta_k \phi_k(X(s)) ds, ~~~t\in [0,1],
\end{equation}
differentiating with respect to $\bs\beta$ we obtain the
the linear differential equations:
\begin{eqnarray}
%\frac{d}{dt}  X_{il}^{\theta_i}(t) &=& X_{il}^{\theta_i}(t)
% e^{\theta_i} \sum_{k=1}^M \beta_k \phi_k'(
%X_{il}(t)) + e^{\theta_i} \sum_{k=1}^M \beta_k \phi_k( X_{il}(t)),
%~~~  X_{il}^{\theta_i}(0) = 0, \label{eq:tilde_X_theta_i}\\
\frac{d}{dt}  X^{\beta_r}(t) &=& X^{\beta_r}(t)
\sum_{k=1}^M \beta_k \phi_k'( X(t)) +
\phi_r(X(t)), ~~~ X^{\beta_r}(0)= 0,\label{eq:tilde_X_beta}
\end{eqnarray}
for $r=1,\ldots,M$, where $X^{\beta_r}(t) := \frac{\partial   X(t)}{\partial \beta_r}$.
The Hessian of $X(\cdot)$ with respect to
$\bs{\beta}$ is given by the matrix $(X^{\beta_r,\beta_{r'}})_{r,r'=1}^M$,
where $X^{\beta_r,\beta_{r'}}(t) :=
\frac{\partial^2}{\partial \beta_r
\partial \beta_{r'}}  X(t)$, which satisfies the system of
ODEs, for $r,r'=1,\cdots,M$:
\begin{eqnarray}\label{eq:tilde_X_beta_Hessian}
&& \frac{d}{dt}  X^{\beta_r,\beta_{r'}}(t) \nonumber\\
&=& \left[ X^{\beta_r,\beta_{r'}}(t)
\sum_{k=1}^M \beta_k \phi_k'( X(t))  +
X^{\beta_r}(t) \phi_{r'}'( X(t))\right. \nonumber\\
&&  +
X^{\beta_{r'}}(t) \phi_{r}'( X(t)) +\left.  X^{\beta_r}(t) X^{\beta_{r'}}(t) \sum_{k=1}^M \beta_k \phi_k^{\prime\prime}( X(t)) \right],
~~{X}^{\beta_r,\beta_{r'}}(0)=0. 
\end{eqnarray}
With $a := X(\delta)$ and $X^a(t)$ denoting $\frac{\partial}{\partial a} X(t)$,
we also have
\begin{equation}\label{eq:tilde_X_a}
\frac{d}{dt} X^a(t) =  g_\beta'(X(t)) X^a(t), ~~X^a(\delta) = 1.
\end{equation}

%\subsubsection*{Expressions when $g$ is positive}

Note that (\ref{eq:tilde_X_beta}), (\ref{eq:tilde_X_beta_Hessian}) and
(\ref{eq:tilde_X_a}) are linear differential equations. If the function
$g_{\bs{\beta}} := \sum_{k=1}^M \beta_k \phi_k$ is positive on the domain then
the gradients of the trajectories can be solved explicitly as follows.
%(see Paul \textit{et al.} (2009) for details):
\begin{eqnarray}
% X_{il}^{\theta_i}(t) &=&  e^{\theta_i} t
%g_{\bs{\beta}}( X_{il}(t)). \label{eq:tilde_X_theta_i_closed}\\
X^{\beta_r}(t) &=& g_{\bs{\beta}}( X(t))
\int_{ x_0}^{ X(t)} \frac{\phi_r(x)}{(g_{\bs{\beta}}(x))^2} dx.
\label{eq:tilde_X_beta_closed}
\end{eqnarray}
\begin{eqnarray}\label{eq:tilde_X_beta_Hessian_alt}
X^{\beta_r,\beta_{r'}}(t) 
&=& g_{\bs{\beta}}(X(t))
\int_0^t \frac{1}{g_{\bs{\beta}}(X(s))}
\left[X^{\beta_r}(s) \phi_{r'}'(X(s)) +  \phi_{r}'(X(s)) X^{\beta_{r'}}(s)\right]ds
\nonumber\\
&& ~+ g_{\bs{\beta}}(X(t))
\int_0^t  \frac{1}{g_{\bs{\beta}}(X(s))} X^{\beta_r}(s)X^{\beta_{r'}}(s)
g_{\bs{\beta}}''(X(s)) ds.  
\end{eqnarray}
and
\begin{equation}\label{eq:X_partial_a}
X^a(t) = \frac{g_{\bs\beta}(X(t))}{g_{\bs\beta}(a)}, ~~~t \in [\delta,1-\delta].
\end{equation}

%\subsubsection*{Application of perturbation theory : rate bounds}\label{subsec:rates}

Now we summarize approximations of various relevant quantities.
The following result on the perturbation of the solution path in an
initial value problem due to a perturbation in the gradient function
is derived from Deuflhard and Bornemann (2002).
%The proof is outlined in Paul \textit{et al.} (2009), Appendix F.
\begin{proposition}\label{prop:perturb}
Consider the initial value problem:
\begin{equation}\label{eq:ODE_general}
x' = f(t,x),~~~x(t_0) = x_0,
\end{equation}
where $x \in \mathbb{R}^d$. On the augmented phase space $\Omega$, say, let the
mappings $f$ and $\delta f$ be continuous and continuously
differentiable with respect to the state variable. Assume that for
$(t_0,x_0) \in \Omega$, the initial value problem
(\ref{eq:ODE_general}), and the perturbed problem
\begin{equation*}\label{eq:ODE_perturbed}
x' = f(t,x) + \delta f(t,x), ~~~x(t_0) = x_0,
\end{equation*}
have the solutions $x$ and $\ol x = x + \delta x$, respectively.
If $f$ is such that
$\parallel f_x(t,\cdot) \parallel_{\infty} \leq \chi(t)$ for a
function $\chi(\cdot)$ bounded on $[t_0,t_1]$, and $\parallel \delta
f(t, \cdot) \parallel_\infty \leq \tau(t)$ for some nonnegative
function $\tau(\cdot)$ on $[t_0,t_1]$, then
\begin{equation*}\label{eq:solution_perturb_norm_bound}
\parallel \delta x(t) \parallel \leq \int_{t_0}^t \exp\left(\int_s^t \chi(u)
du\right) \tau(s) ds, ~~~~\mbox{for all}~t \in [t_0,t_1].
\end{equation*}
\end{proposition}
We use the above result to compute bounds for the trajectories and
their derivatives corresponding to the different values of the parameter
$\bs\beta$ in a neighborhood of the point $\bs\beta^*$. In order to keep the
exposition simple, we assume that $g_{\bs\beta}(x) =
g_{\bs\beta}(x_{1,M})$ for $x > x_{1,M}$ and $g_{\bs\beta}(x) =
g_{\bs\beta}(x_{0,M})$ for $x < x_{0,M}$ with a differentiability requirement
at the points $x_{0,M}$ and $x_{1,M}$.

Our aim is to show that the range of the trajectories
$X(t;\bs\beta,x_{0,\delta})$ is contained in the set $D_0 = [x_{0,M},x_{1,M}]$,
for all $t \in [\delta,1-\delta]$ and for all $\bs\beta \in {\cal B}(\alpha_n)
:= \{\bs\beta:
\parallel \bs\beta - \bs\beta^*\parallel \leq \alpha_n\}$.  Let $\gamma_n =
\max\{\sup_{x\in D} |g_{\bs\beta^*}(x) - g(x)|, \sup_{x\in D_0}
|g_{\bs\beta^*}(x) - g_{\bs\beta}(x)|\}$. Then $\gamma_n = O(M^{-p}) +
O(\alpha_n M^{1/2})$. Also, let $\xi_n = \max_{j=0,1}|\widehat x_j -
x_{j,\delta}|$. As in the proof of Proposition \ref{prop:kappa}, we can easily
show that $[x_{0,\delta},x_{1,\delta}] \subset [x_{0,M},x_{1,M}]$ for
sufficiently large $M$ . On the other hand, by using the perturbation bound given by Proposition \ref{prop:perturb}
progressively over small subintervals of the interval $[\delta,1-\delta]$, it
can be shown that
\begin{equation*}
\sup_{\bs\beta \in {\cal B}(\alpha_n)} \sup_{t \in
[\delta,1-\delta]}|X(t;\bs\beta,x_{0,\delta}) - X_g(t;x_{0,\delta})| \leq C_1
\gamma_n + C_2 \xi_n,
\end{equation*}
for appropriate positive constants $C_1,C_2$ that depend on the value of $g$
and $g'$ on the interval $[x_0,x_1]$. Now, using
%the bound on $\xi_n$ stated in
Lemma \ref{lem:consistency_x_j_hat}, the condition on $\alpha_n$ as given in
Theorem \ref{thm:consistency}, and the definitions of $\widehat x_j$,
$x_{j,\delta}$ and $x_{j,M}$, for $j=0,1$, we conclude that for large enough
$M$, the range of $X(t;\bs\beta,x_{0,\delta})$ is contained in $D_0$ for all $t
\in [\delta,1-\delta]$ and for all $\bs\beta \in {\cal B}(\alpha_n)$. The
scenario is depicted in Figure \ref{fig:trajectory_envelop},
where the dashed curves indicate the envelop of the trajectories
$X(t;\bs\beta,x_{0,\delta})$, while the solid curve indicates the trajectory
$X_g(t;x_{0,\delta})$.

%\textcolor{green}{[ADD THE FIGURE AND THE CORRESPONDING PERTURBATION THEORETIC
%ARGUMENT FOR $[x_{0,\delta},x_{1,\delta}] \subset [x_{0,M},x_{1,M}]$].}

Next, we provide bounds for trajectories and their derivatives. In the following,
$\parallel \cdot \parallel_\infty$ is used to denote the $\sup$-norm over $D_0 =
[x_{0,M},x_{1,M}]$. First, by {\bf A2} we have the following:
\begin{equation}\label{eq:g_beta_estimates}
\parallel g_{\bs{\beta}}^{(j)} - g_{\bs{\beta}^*}^{(j)} \parallel_\infty
= O(\parallel \bs\beta - \bs\beta^*\parallel M^{j+1/2})  ~~~j=0,1,2,
\end{equation}
where $g^{(j)}$ and $g_{\bs{\beta}^*}^{(j)}$ denote the $j$-th derivative of
$g$ and $g_{\bs{\beta}^*}$, respectively. Next, again from {\bf A2}, for $M$
large enough, solutions $\{X(t;\bs{\beta}):t \in [\delta,1-\delta]\}$ exist for
all $\bs{\beta}$ such that $\parallel \bs{\beta} - \bs{\beta}^*\parallel \leq
\alpha_n$. This also implies that the solutions $X^{\beta_r}(\cdot;\bs{\beta})$
and $X^{\beta_r,\beta_{r'}}(\cdot;\bs{\beta})$ exist on $[\delta,1-\delta]$ for
all $\bs{\beta}$ such that $\parallel \bs{\beta} - \bs{\beta}^*\parallel \leq
\alpha_n$, since they follow linear differential equations where the
coefficient functions depend on $X(t;\bs{\beta})$. Moreover, by
\textit{Gronwall's lemma} (Deuflhard and Bornemann, 2002),
(\ref{eq:g_beta_estimates}) and the fact that $\parallel
g_{\bs{\beta}^*}^{(j)}\parallel_\infty = O(1)$ for $j=0,1,2$ (again by {\bf
A2}).

Hence, if $\parallel \bs\beta - \bs\beta^*\parallel M^{3/2} = o(1)$, then using Proposition \ref{prop:perturb},
the fact that $\parallel g_{\bs{\beta}^*}^{(j)}\parallel_\infty = O(1)$ for
$j=0,1,2$, and the expressions for the ODEs for the partial derivatives,
%after some algebra
we obtain (almost surely):
\begin{equation}\label{eq:X_path_bias}
\parallel X(\cdot;\bs{\beta}^*) - X_g(\cdot)\parallel_\infty
= O(M^{-p}).
\end{equation}
The same technique can be used to prove the following:
\begin{eqnarray}
\parallel X(\cdot;\bs{\beta}) - X(\cdot;\bs{\beta}^*)\parallel_\infty
&=& O(\parallel \bs\beta - \bs\beta^*\parallel M^{1/2})   \label{eq:X_diff_bound}\\
%\parallel X_{il}^{\theta_i}(\cdot;\bs{\theta},\bs{\beta}) - X_{il}^{\theta_i}(\cdot;\bs{\theta}^*,\bs{\beta}^*)
%\parallel_\infty
%&=& O(\alpha_N M^{3/2}) \label{eq:X_diff_theta_bound}\\
\max_{1\leq r\leq M} \parallel X^{\beta_r}(\cdot;\bs{\beta}^*) \parallel_\infty
&=& O(M^{-1/2}) \label{eq:X_beta_bound}\\
\max_{1\leq r\leq M} \parallel X^{\beta_r}(\cdot;\bs{\beta}) - X^{\beta_r}(\cdot;\bs{\beta}^*)
\parallel_\infty
&=& O(\parallel \bs\beta - \bs\beta^*\parallel M) \label{eq:X_diff_beta_bound}\\
%\parallel X_{il}^{\theta_i,\theta_i}(\cdot;\bs{\theta},\bs{\beta}) - X_{il}^{\theta_i,\theta_i}(\cdot;\bs{\theta}^*,\bs{\beta}^*)
%\parallel_\infty
%&=& O(\alpha_N M^{5/2}) \label{eq:X_diff_theta_theta_bound}\\
%\max_{1\leq r\leq M} \parallel X_{il}^{\theta_i,\beta_r}(\cdot;\bs{\theta}^*,\bs{\beta}^*) \parallel_\infty
%&=& O(M^{1/2}) \label{eq:X_theta_beta_bound}\\
%\max_{1\leq r\leq M} \parallel X_{il}^{\theta_i,\beta_r}(\cdot;\bs{\theta},\bs{\beta})
%- X_{il}^{\theta_i,\beta_r}(\cdot;\bs{\theta}^*,\bs{\beta}^*) \parallel_\infty
%&=& O(\alpha_N M^2) \label{eq:X_diff_theta_beta_bound}\\
\max_{1\leq r,r'\leq M} \parallel X^{\beta_r,\beta_{r'}}(\cdot;\bs{\beta}^*)
\parallel_\infty &=& O(1) \label{eq:X_beta_beta_bound}\\
\max_{1\leq r,r'\leq M} \parallel X^{\beta_r,\beta_{r'}}(\cdot;\bs{\beta}) -
X^{\beta_r,\beta_{r'}}(\cdot;\bs{\beta}^*) \parallel_\infty
&=& O(\parallel \bs\beta - \bs\beta^*\parallel M^{3/2}) \label{eq:X_diff_beta_beta_bound}
\end{eqnarray}
whenever $\parallel \bs{\beta} - \bs{\beta}^*\parallel M^{3/2} = o(1)$.

To illustrate the key arguments, we prove (\ref{eq:X_beta_bound}) and
(\ref{eq:X_diff_beta_bound}). First, (\ref{eq:X_beta_bound}) follows by
(\ref{eq:tilde_X_beta_closed}), and the fact that $\parallel \phi_{r}
\parallel_\infty = O(M^{1/2})$ and is supported on an interval of length
$O(M^{-1})$. In fact it holds for all $\bs\beta$ such that
$\parallel \bs{\beta} - \bs{\beta}^*\parallel M^{3/2} = o(1)$. Next, note that the
function $\phi_{r}$  is Lipschitz with Lipschitz constant $O(M^{3/2})$ and is
supported on an interval of length $O(M^{-1})$. Since (\ref{eq:tilde_X_beta})
is a linear differential equation, using Proposition \ref{prop:perturb} with
$\delta f(t,x)$ given by
\begin{eqnarray*}
&& x\left[g_{\bs{\beta}}'(X(t;\bs{\theta},
\bs{\beta})) - g_{\bs{\beta}^*}'(X(t;\bs{\beta}^*))
\right]
+ \phi_{r}(X(t;\bs{\beta})) - \phi_{r}(X(t;\bs{\beta}^*))
\end{eqnarray*}
we obtain (\ref{eq:X_diff_beta_bound}) by using (\ref{eq:X_diff_bound}) and the following facts:
$\sup_{t\in [\delta,1-\delta]} |X^{\beta_r}(t;\bs\beta)| = O(M^{-1/2})$
for all $\bs{\beta} \in \Omega(\alpha_n)$;
$\parallel g_{\bs{\beta}}''\parallel_\infty = O(\alpha_n M^{5/2})$;
$\parallel g_{\bs{\beta}}' - g_{\bs{\beta}^*}'
\parallel_\infty = O(\alpha_n M^{3/2})$; and $\alpha_n M^{3/2} = o(1)$.

\subsection*{Appendix B : Verification of (vi) of {\bf A2} for B-spline basis}\label{subsec:spline_approx}

In this subsection, we verify that the condition (vi) of {\bf A2} is satisfied if $\{\phi_{k,M}\}_{k=1}^M$ is a
normalized B-spline basis with equally spaced knots on $[x_{0,M},x_{1,M}]$ and of order $d\geq \max\{3,p-1\}$.
In particular, we show that the rate of approximation of $X(t)$ by
$X(t;\bs\beta^*)$  with a carefully chosen $\bs\beta = \bs\beta^*$ satisfies
the requirement that $\sup_{t \in [\delta,1-\delta]} |X(t) - X(t;\bs\beta^*)| = O(M^{-(p+1)})$ and
the conditions $\sup_{x \in [x_{0,M},x_{1,M}]} |g^{(j)}(x) - g_{\bs\beta^*}^{(j)}(x)| = O(M^{-p+j})$ for
$j=0,1,2$. The result is proved through the following lemmas proved in SM.

\begin{lemma}\label{lemma:path_bound_refined}
Suppose that $\{\phi_{k,M}\}_{k=1}^M$ has combined support $[x_{0,\delta}, x_{1,\delta}] = [X(\delta),X(1-\delta)]$ and
satisfies (ii)--(v) of {\bf A2} and $\bs\beta^*$ furthermore has the property that
\begin{equation}\label{eq:g_inverse_intergal_bound}
\sup_{x \in [x_{0,\delta},x_{1,\delta}]} \left|\int_{x_{0,\delta}}^{x} \frac{g(u) - g_{\bs\beta^*}(u)}{g(u)} du  \right| = a_M
\end{equation}
such that $c_0 M^{-(p+1)} \leq a_M \ll M^{-p-\epsilon}$, uniformly in $M$, for some $\epsilon \in (0,1]$ and some
$c_0 > 0$.  Then, if $X(\delta;\bs\beta^*) = X(\delta)$,
there exists $C > 0$ such that
\begin{equation}\label{eq:X_diff_bound_refined}
\sup_{t \in [\delta,1-\delta]} |X(t)-X(t;\bs\beta^*)|  \leq C a_M.
\end{equation}
\end{lemma}

\begin{lemma}\label{lemma:integral_g_spline_approx}
Suppose that {\bf A1} holds with $p\geq 2$. Let $\{\phi_{k,M}\}_{k=1}^M$ denotes the normalized B-spline basis
of order $\geq (p-1)$ with equally spaced knots on the interval $[x_{0,M},x_{1,M}]$.
Then there exists a $\bs\beta^* \in \mathbb{R}^M$ such that
$g_{\bs\beta^*} = \sum_{k=1}^M \beta_k^*\phi_{k,M}$ satisfies
\begin{equation}\label{eq:integral_g_spline_approx}
\sup_{x \in [x_{0,\delta},x_{1,\delta}]} \left|\int_{x_{0,\delta}}^{x}
\frac{g(u) - g_{\bs\beta^*}(u)}{g(u)} du  \right| = O(M^{-(p+1)}).
\end{equation}
\end{lemma}

\begin{figure}[th]
\begin{center}
\includegraphics[width=3in, height=2.5in, bb = 100 50 600 600]{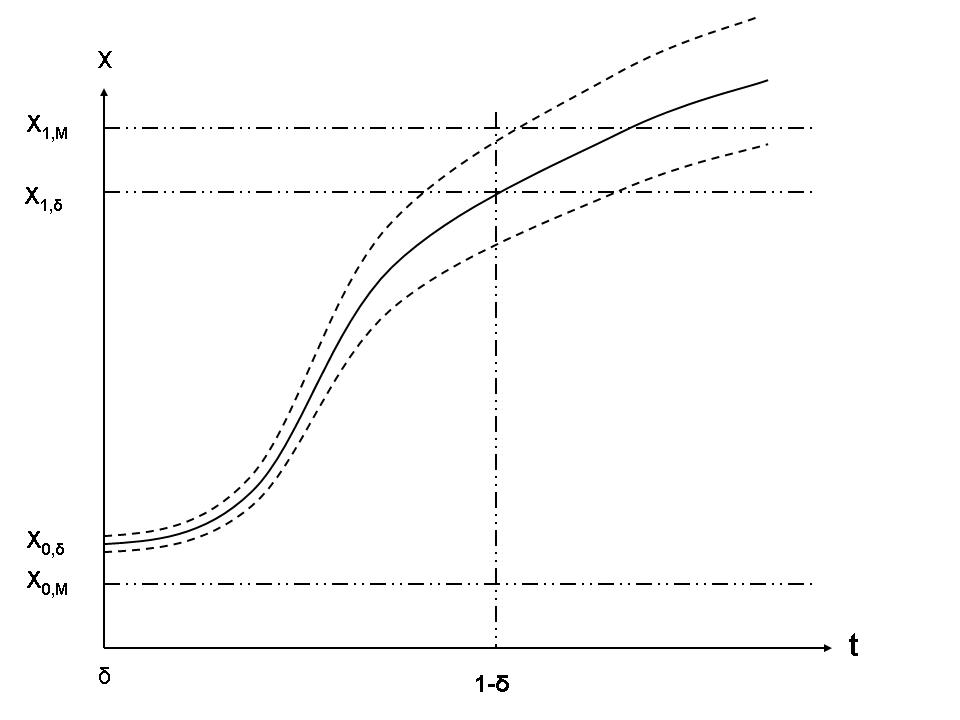}
\caption{Schematic diagram of the trajectory $X_g(t;x_{0,\delta})$ (solid
curve) and the envelop of trajectories $X(t;\bs\beta,x_{0,\delta})$ (boundaries
indicated by dashed curves).
\label{fig:trajectory_envelop}}
\end{center}
\end{figure}

\subsection*{Appendix C :  Proof of Lemma \ref{lemma:Gamma_beta_beta_star_bound}}\label{sec:Appendix_C}

By a Taylor expansion we have, for $j=1,\ldots,n$,
\begin{eqnarray*}
X(T_j;\bs\beta) - X(T_j;\bs\beta^*) &=&  X^{\bs\beta}(T_j;\bs\beta^*)^T(\bs\beta - \bs\beta^*) \\
&& + (X^{\bs\beta}(T_j;\tilde{\bs\beta}(T_j)) - X^{\bs\beta}(T_j;\bs\beta^*))^T(\bs\beta - \bs\beta^*),
\end{eqnarray*}
where $\parallel \tilde{\bs\beta}(T_j) - \bs\beta^*\parallel \leq \parallel \bs\beta - \bs\beta^*\parallel$ for
all $j$. From this, it follows that, for all $\bs\beta \in {\cal A}_M(\alpha_n,\bar\alpha_n)$,
\begin{eqnarray}\label{eq:Gamma_bound_preliminary}
\Gamma_n(\bs\beta,\bs\beta^*) 
&\geq& \frac{3}{4} (\bs\beta - \bs\beta^*)^T \left[\frac{1}{n} \sum_{j=1}^n X^{\bs\beta}(T_j;\bs\beta^*)
X^{\bs\beta}(T_j;\bs\beta^*)^T \mathbf{1}_{[\delta,1-\delta]}(T_j)\right] (\bs\beta - \bs\beta^*) \nonumber\\
&& - 3 \parallel \bs\beta - \bs\beta^*\parallel^2 \frac{1}{n} \sum_{j=1}^n \parallel
X^{\bs\beta}(T_j;\tilde{\bs\beta}(T_j))- X^{\bs\beta}(T_j;\bs\beta^*)\parallel^2 \mathbf{1}_{[\delta,1-\delta]}(T_j),
\end{eqnarray}
where we have used $|2ab| \leq a^2/4 + 4b^2$. Using Proposition \ref{prop:kappa} and Lemma \ref{lem:G_gamma_quad_bound}
(stated below) we conclude, given $\eta > 0$, there exists $C_{10}(\eta) > 0$ such that,
\begin{eqnarray*}
&& (\bs\beta - \bs\beta^*)^T \left[\frac{1}{n} \sum_{j=1}^n X^{\bs\beta}(T_j;\bs\beta^*)
X^{\bs\beta}(T_j;\bs\beta^*)^T \mathbf{1}_{[\delta,1-\delta]}(T_j) \right] (\bs\beta - \bs\beta^*) \\
&\geq& C_{10}(\eta) \frac{1}{M^2} \parallel \bs\beta - \bs\beta^*\parallel^2
\end{eqnarray*}
for all $\bs\beta \in {\cal A}_M(\alpha_n,\bar\alpha_n)$,
with probability at least $1-n^{-\eta}$. Now, another application of the Mean Value Theorem
yields that for $T_j \in [\delta,1-\delta]$,
\begin{eqnarray*}
&& \parallel X^{\bs\beta}(T_j;\tilde{\bs\beta}(T_j))- X^{\bs\beta}(T_j;\bs\beta^*)\parallel^2 \\
&\leq& \parallel \tilde{\bs\beta}(T_j) - \bs\beta^*\parallel^2 \parallel X^{\bs\beta\bs\beta^T}(T_j;\bs\beta^*)\parallel_F^2
\\
&&  + \parallel \tilde{\bs\beta}(T_j) - \bs\beta^*\parallel^2 \sum_{1\leq k,k' \leq M}|X^{\beta_k,\beta_{k'}}(T_j;\bar{\bs\beta}^{k}(T_j)) -  X^{\beta_k,\beta_{k'}}(T_j;\bs\beta^*)|^2,
\end{eqnarray*}
where $\parallel \cdot \parallel_F$ denotes the Frobenius norm, and
$\parallel \bar{\bs\beta}^{k}(T_j) - \bs\beta^*\parallel \leq \parallel \tilde{\bs\beta}(T_j) - \bs\beta^*\parallel$
for all $1\leq k \leq M$ and $1\leq j\leq n$. Now, using (\ref{eq:X_beta_beta_bound}) and (\ref{eq:X_diff_beta_beta_bound}),
and combining the last three displays, we get (\ref{eq:Gamma_beta_beta_star_bound}).

\begin{lemma}\label{lem:G_gamma_quad_bound}
Suppose that {\bf A1}--{\bf A4} hold. Let
$$
\bar{G}_{*n} := \frac{1}{F_T(1-\delta) - F_T(\delta)} \frac{1}{n}
\sum_{j=1}^n X^{\bs\beta}(T_j;\bs\beta^*)(X^{\bs\beta}(T_j;\bs\beta^*))^T
\mathbf{1}_{[\delta,1-\delta]}(T_j).
$$
Then, given $\eta > 0$, there exists  constants $c_1'(\eta), c_2'(\eta) > 0$  such that, with probability $1-n^{-\eta}$,
uniformly in $\bs{\gamma} \in \mathbb{S}^{M-1}$,
\begin{equation}\label{eq:G_gamma_quad_bound}
\bs{\gamma}^T \bar{G}_{*n} \bs{\gamma} \geq
\bs{\gamma}^T G_* \bs{\gamma}  - c_1'(\eta)\sqrt{\bs{\gamma}^T G_* \bs{\gamma}}\sqrt{\frac{M\log n}{n}} \geq
c_2'(\eta) M^{-2}.
\end{equation}
\end{lemma}

\subsubsection*{Proof of Lemma \ref{lem:G_gamma_quad_bound}}
Let $\mathbf{v}_j =
X^{\bs\beta}(T_j;\bs\beta^*)$. Define $D(\bs{\gamma}) = \bs{\gamma}^T (\bar G_{*n}  - G_*) \bs{\gamma}$.
Notice that
$$
\frac{1}{(F_T(1-\delta)-F_T(\delta))}\mathbb{E}_{T}[\mathbf{v}_j\mathbf{v}_j^T \mathbf{1}_{[\delta, 1-\delta]}(T_j)]
= \mathbb{E}_{\tilde T}[\mathbf{v}_j\mathbf{v}_j^T]
= G_*,
$$
where the first expectation is with respect to the distribution of $T_1$ and the second with respect to that of $\tilde T_1$.
Hence, we can write
$D(\bs\gamma) = n^{-1}\sum_{j=1}^{n} u_j(\bs\gamma)$
where
$$
u_j(\bs\gamma) = \bs\gamma^T \left(\mathbf{v}_j\mathbf{v}_j^T \frac{\mathbf{1}_{[\delta, 1-\delta]}(T_j)}{F_T(1-\delta)-F_T(\delta)}
- \mathbb{E}_{\tilde T}[\mathbf{v}_j\mathbf{v}_j^T]\right)\bs\gamma.
$$
%and
%$w_{ilj}(\bs\gamma) = \bs\gamma^T (\mathbb{E}[(\nabla_{i}
%X_{ilj}\nabla_{i} X_{ilj}^T)|a_{il}]  - \mathbb{E}[\nabla_{i}
%X_{ilj}\nabla_{i} X_{ilj}^T)])\bs\gamma$, where, for notational
%simplicity,
%$\nabla_i X_{ilj} := ( (X_{il}^{\bs{\beta}}(T_{i,j}))^T, (X_{il}^{\theta_i}(T_{i,j}) \mathbf{e}_{i-1})^T)^T$,
%$i=1,\ldots,n$.
Note that, the random variables $u_j(\bs\gamma)$ have zero
conditional mean, are uniformly bounded, and are independent.
%and the variables $Z_j(\bs\gamma) := \sum_{l=1}^{N_i} u_{ilj}(\bs\gamma)$
%are independent.
%Similarly, the random variables $\{w_{ilj}\}_{i,j}$
%have zero mean, are uniformly bounded, and the variables
%$\sum_{j=1}^{m_i} w_{ilj}(\bs\gamma)$ are independent. Indeed, for
%each fixed $(i,l)$, the variables $\{w_{ilj}\}_{j=1}^{m_i}$ are
%identical.
% since $T_{i,j}$ are i.i.d.
Moreover, the functions $u_j(\bs\gamma)$
are differentiable functions of $\bs\gamma$.
%Define ${\cal G}_*(\bs{a}) :=
%\mathbb{E}({\cal G}_*(\bs\theta^*,\bs\beta^*)|\bs{a})$.
Then, since by (\ref{eq:X_beta_bound}), $u_j(\bs\gamma)$'s are uniformly bounded by some
$K_1 > 0$,
\begin{eqnarray*}
\mbox{Var}\left(\sum_{j=1}^{n} u_j(\bs\gamma)\right) ~=~
\sum_{j=1}^{n} \mathbb{E}[(u_j(\bs\gamma))^2]
&\leq& K_1 \sum_{j=1}^{n}\mathbb{E}|u_j(\bs\gamma)| ~\leq~
2K_1 n \bs\gamma^T G_* \bs\gamma .
\end{eqnarray*}
%\begin{eqnarray*}
%\mbox{Var}(\sum_{i=1}^n\sum_{j=1}^{m_i} Z_{ij}(\bs\gamma) |\bs{a}) &=&
%\sum_{i=1}^n\sum_{j=1}^{m_i} \mathbb{E}[(Z_{ij}(\bs\gamma))^2|\bs{a}]
%~\leq~ \sum_{i=1}^n\sum_{j=1}^{m_i} N_i \sum_{l=1}^{N_i} \mathbb{E}[u_{ilj}^2(\bs\gamma)|a_{il}] \nonumber\\
%&\leq& K_2\ol{N} \sum_{i=1}^n\sum_{l=1}^{N_i}\sum_{j=1}^{m_i}\mathbb{E}[|u_{ilj}(\bs\gamma)||a_{il}] ~\leq~
%2K_2\ol{N} \bs\gamma^T{\cal
%G}_*(\bs{a}) \bs\gamma .
%\end{eqnarray*}
%In the above, second inequality uses  $(\sum_{i=1}^N x_i)^2 \leq
%N\sum_{i=1}^N x_i^2$, and the last follows from fact that
%$u_{ilj}(\bs\gamma)$ is a difference of two nonnegative quantities,
%the second one being the conditional expectation of the first one
%given $\bs{a}$.
Thus, by Bernstein's inequality, for every $v
> 0$ and $\bs\gamma \in \mathbb{S}^{M-1}$,
\begin{equation*}
\mathbb{P}\left(|\sum_{j=1}^{n}
u_j(\bs\gamma)| > v\right) \leq 2\exp\left(-\frac{v^2/2}{2K_1 n \bs\gamma^T
G_* \bs\gamma + K_1 v/3}\right).
\end{equation*}
On the other hand, by (\ref{eq:G_star_condition}), $\bs\gamma^T G_* \bs\gamma \geq c M^{-2}$ for some $c > 0$.
By this, and the condition that $M^3 = o(n/\log n)$, it is easy to see that $\sqrt{\bs\gamma^T G_*  \bs\gamma}
\gg \sqrt{M \log n/n}$.
Thus, using  an entropy argument as in the proof of (\ref{eq:Z_beta_sup_bound}), we
conclude that given $\eta > 0$ there exists
$c_1'(\eta) > 0$ such that
%on the set $\{\bs{a} | \bs\gamma^T{\cal G}_*(\bs{a}) \bs\gamma \geq C_2(\delta) \ol{N} M
%\log(\ol{N}\ol{m})\}$,
\begin{equation}\label{eq:sum_u_bound}
\mathbb{P}\left(\sup_{\bs\gamma \in \mathbb{S}^{M-1}} \frac{|n^{-1}\sum_{j=1}^{n}
u_j(\bs\gamma)|}{\sqrt{\bs\gamma^T G_* \bs\gamma}} \leq c_1'(\eta) \sqrt{\frac{M
\log n}{n}} \right) > 1 - n^{-\eta}.
\end{equation}
Recalling the definition of $D(\bs\gamma)$,  and again using the fact that
$\bs\gamma^T G_* \bs\gamma \geq c M^{-2}$ and $M^3 = o(n/\log n)$, (\ref{eq:G_gamma_quad_bound})
follows from (\ref{eq:sum_u_bound}).

\bibliographystyle{plain}

\clearpage
\newpage

\setcounter{figure}{0}
\renewcommand{\thefigure}{\arabic{figure}}

\begin{figure}[th]
\begin{center}
\begin{tabular}{llllllll}
\includegraphics[width=1.3in, height=1.5in, bb = 0 0 600 600]{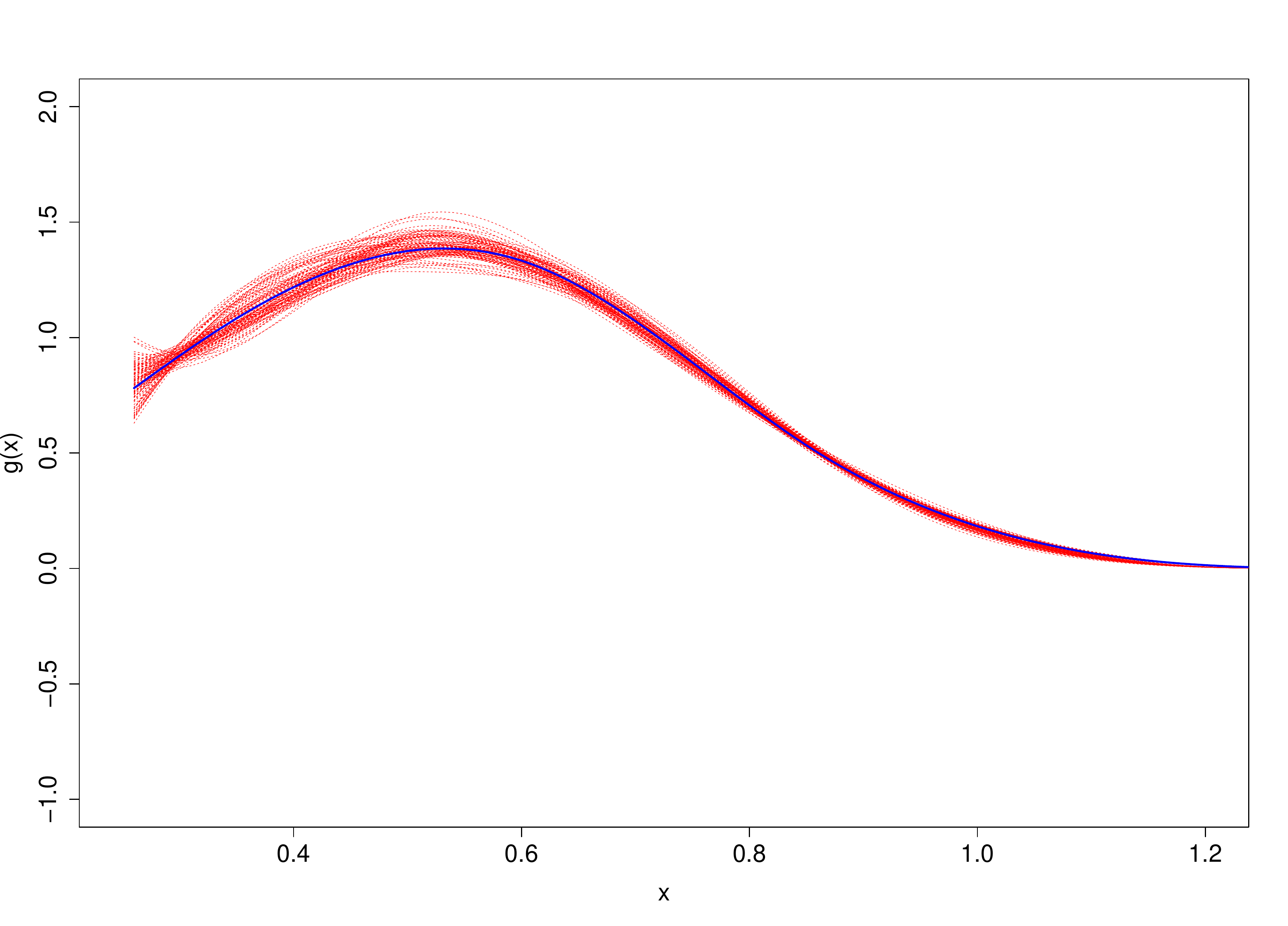} &&&&&&
\includegraphics[width=1.3in, height=1.5in, bb = 0 0 600 600]{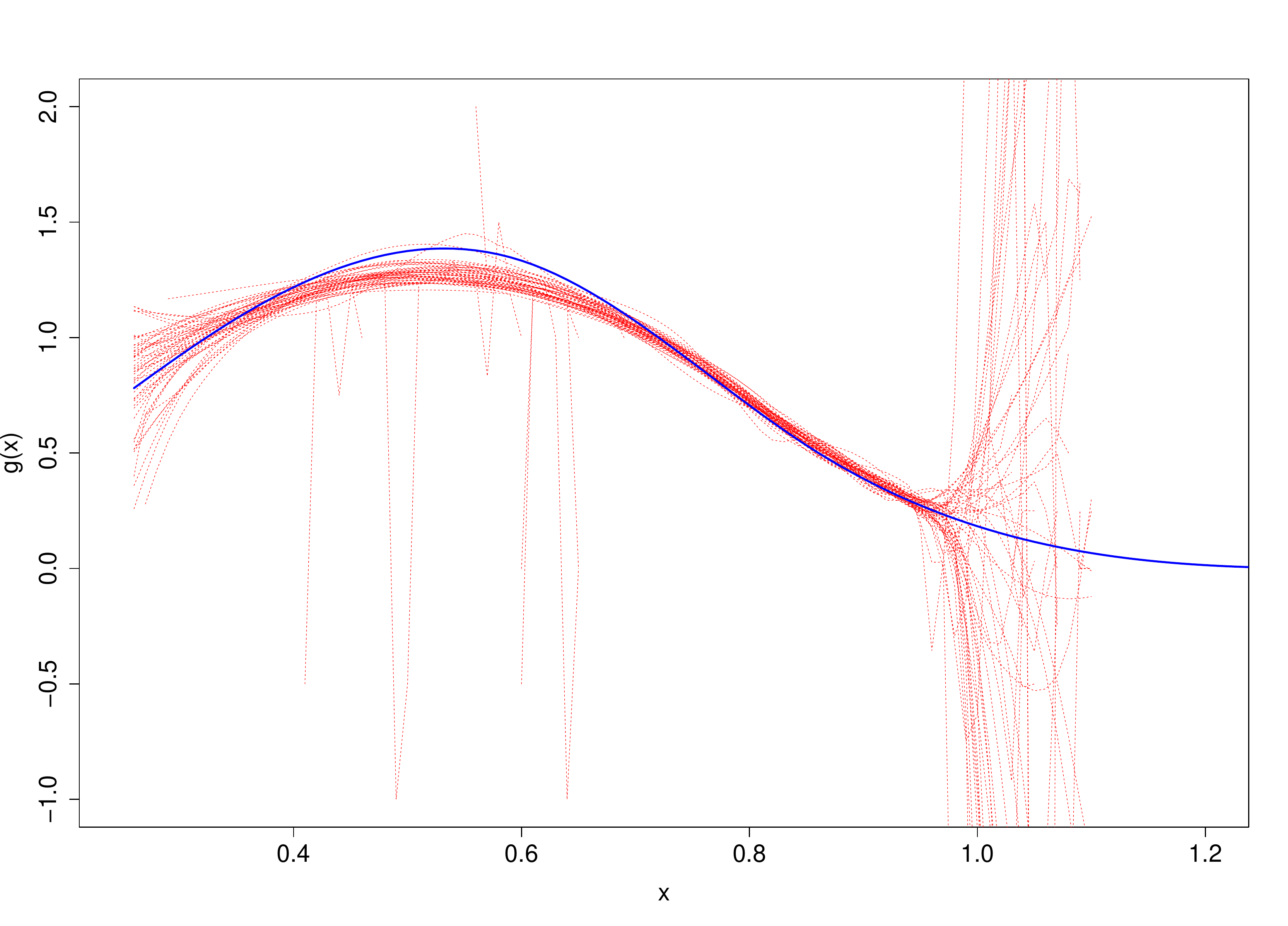}\\
\end{tabular}
\caption{Simulation: Estimated gradient functions (red curves) overlayed on the
true gradient function (blue curve). Left panel: proposed estimator; Right
panel: two-stage estimator.
\label{fig:simu}}
\end{center}
\end{figure}

\begin{figure}[th]
\begin{center}
\includegraphics[width=2in, height=2in, bb = 100 0 700 700]{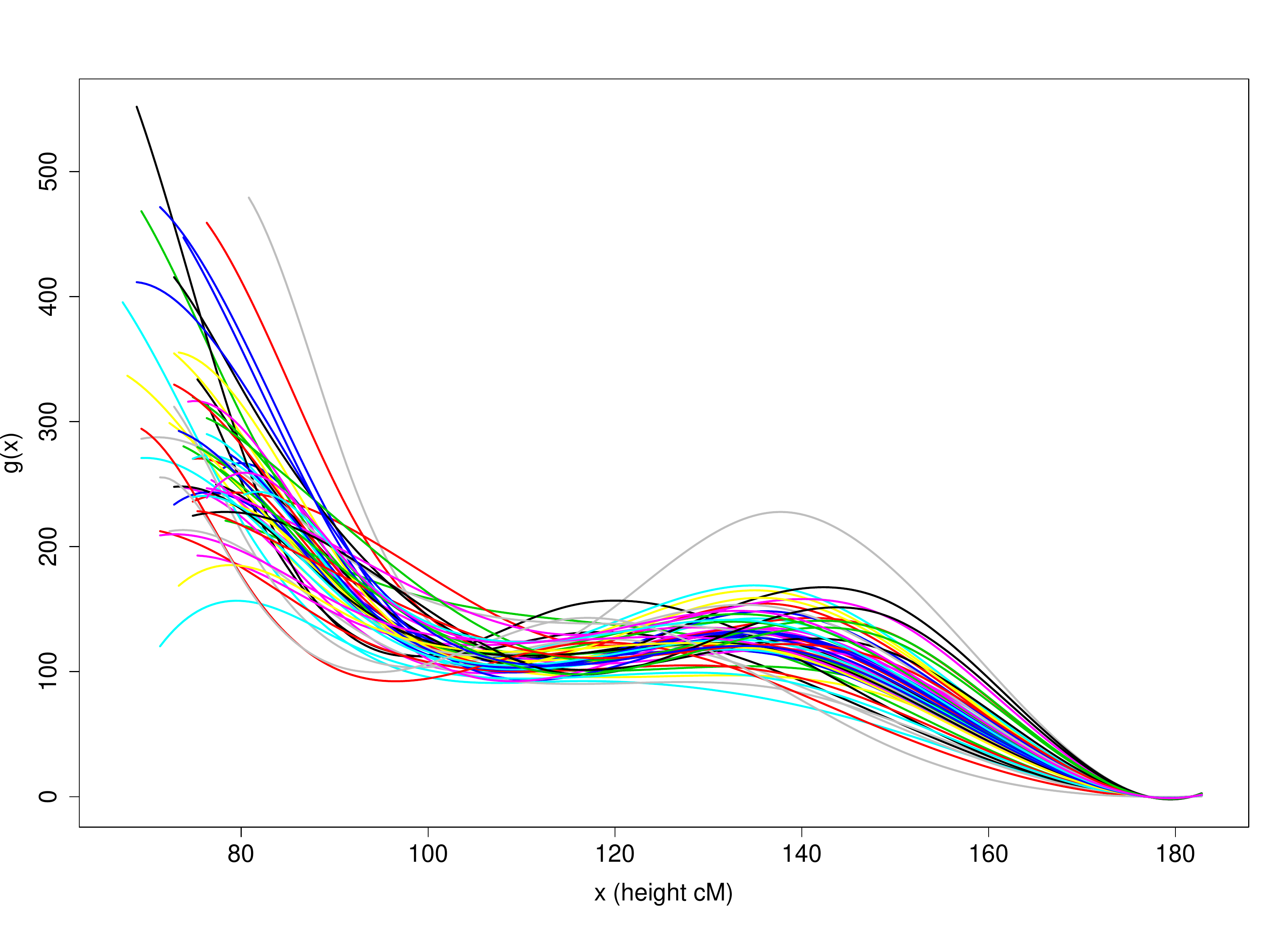}
\caption{Berkeley Growth Data: fitted gradient functions for $54$ female
subjects.
\label{fig:fit_gradient}}
\end{center}
\end{figure}

\clearpage
\newpage

%%%%%%%%%%%%%%%%%%%%%%%%%%%%%%%%%%%%%%%%%%%%%%%%%%%%%%%
%%%%%%%%%% Supplementary Material %%%%%%%%%%%%%%%%%%%%%
%%%%%%%%%%%%%%%%%%%%%%%%%%%%%%%%%%%%%%%%%%%%%%%%%%%%%%%

%\begin{center}
%\Large{\bf Supplementary Materials for ``Consistency in a semiparametric
%autonomous dynamical system'' }
%\end{center}

\setcounter{page}{1}

\setcounter{section}{0}
\renewcommand{\thesection}{S\arabic{section}}
\setcounter{equation}{0}
\renewcommand{\theequation}{S.\arabic{equation}}
\setcounter{figure}{0}
\renewcommand{\thefigure}{S.\arabic{figure}}
\setcounter{proposition}{0}
\renewcommand{\theproposition}{S.\arabic{proposition}}
\setcounter{lemma}{0}
\renewcommand{\thelemma}{S.\arabic{lemma}}
\setcounter{corollary}{0}
\renewcommand{\thecorollary}{S.\arabic{corollary}}

\begin{center}
\Large{\bf Supplementary Material : ``Nonparametric estimation of dynamics of monotone trajectories'' }
\end{center}

\vskip.1in

\section{Proof of Lemmas \ref{lemma:path_bound_refined} and \ref{lemma:integral_g_spline_approx}}

\subsubsection*{Proof of Lemma \ref{lemma:path_bound_refined}}

First, we write (since $X(\delta;\bs\beta^*) = X(\delta) = x_{0,\delta}$), for $t \in [\delta,1-\delta]$,
\begin{eqnarray}\label{eq:X_X_beta_diff_expand}
X(t) - X(t;\bs\beta^*) 
&=& \int_\delta^t (g(X(s)) - g_{\bs\beta^*}(X(s;\bs\beta^*)))ds \nonumber\\
&=& \int_{\delta}^t (g(X(s)) - g_{\bs\beta^*}(X(s)))ds + \int_{\delta}^t (g_{\bs\beta^*}(X(s)) - g_{\bs\beta^*}(X(s;\bs\beta^*)))ds \nonumber\\
&=& \int_{X(\delta)}^{X(t)} (g(u) - g_{\bs\beta^*}(u)) \frac{du}{g(u)} + \int_{\delta}^t (g_{\bs\beta^*}(X(s)) - g_{\bs\beta^*}(X(s;\bs\beta^*)))ds.
\end{eqnarray}
In the last step we have used $X'(s) = g(X(s))$. Since $\sup_{x \in [x_{0,\delta},x_{1,\delta}]}
|g(x) - g_{\bs\beta^*}(x)| = O(M^{-p})$,  from (A.10) we already have $\sup_{t \in [\delta,1-\delta]}|X(t) - X(t;\bs\beta^*)| \leq C_1 M^{-p}$
for some $C_1 > 0$. Also, $\sup_{x \in [x_{0,\delta},x_{1,\delta}]\cap \Pi^c} \parallel g_{\bs\beta^*}'(x) \parallel \leq C_2$ for some $C_2 > 0$,
where $\Pi$ denotes the set of knots for the B-spline functions. This implies, by Mean Value Theorem,
that $g_{\bs\beta^*}$ is a Lipschitz function with Lipschitz constant bounded by $C_2$.
Thus, for all $t \in [\delta,1-\delta]$,
\begin{eqnarray}\label{eq:diff_g_beta_X_beta}
\int_{\delta}^t (g_{\bs\beta^*}(X(s)) - g_{\bs\beta^*}(X(s;\bs\beta^*)))ds
&\hskip-.1in\leq& \hskip-.1in C_2  \int_{\delta}^t |X(s) - X(s;\bs\beta^*)| ds \nonumber\\
&\hskip-.1in\leq& \hskip-.1in C_1 C_2 M^{-p} (t-\delta). 
\end{eqnarray}
Combining (\ref{eq:g_inverse_intergal_bound}) and (\ref{eq:diff_g_beta_X_beta}), we have
\begin{equation*}
|X(t) - X(t;\bs\beta^*)| \leq a_M + C_1 C_2 M^{-p}(t-\delta), \qquad t \in [\delta,1-\delta].
\end{equation*}
Substituting this again in the last line of (\ref{eq:diff_g_beta_X_beta}), from (\ref{eq:X_X_beta_diff_expand}), we obtain
\begin{eqnarray*}
|X(t) - X(t;\bs\beta^*)| &\leq& a_M + C_2 \int_\delta^t (a_M + C_1 C_2 M^{-p}(s-\delta))ds \\
&=& a_M + C_2 a_M(t-\delta) + C_1 C_2^2 M^{-p} \frac{(t-\delta)^2}{2!}
\end{eqnarray*}
for $t \in [\delta,1-\delta]$. By induction, it follows that for all $J \geq 1$,
\begin{equation*}\label{eq:X_X_beta_diff_final}
|X(t) - X(t;\bs\beta^*)| \leq a_M \sum_{j=0}^{J} C_2^j \frac{(t-\delta)^j}{j!} + C_1 C_2^{J+1} \frac{(t-\delta)^{J+1}}{(J+1)!}, ~~t \in [\delta,1-\delta].
\end{equation*}
Since $c_0 M^{-(p+1)} < a_M \ll M^{-p-\epsilon}$, we obtain (\ref{eq:X_diff_bound_refined})
by choosing $J$ sufficiently large and recalling the
expansion of $e^{C_2(t-\delta)}$, whereby we can take $C = 2e^{C_2}$.

\subsubsection*{Proof of Lemma \ref{lemma:integral_g_spline_approx}}

Define $G(x) = \int_{x_{0,M}}^x g(u) du$ for $x \in [x_{0,M},x_{1,M}]$.
It is well known (cf. de Boor, 1978, ch. XII) that, for every $d \geq p$,  there exists a spline $S_d(x)$
of order $d \geq p$ with equally spaced knots
with spacing $O(M^{-1})$ on the interval $[x_{0,M},x_{1,M}]$ such that
\begin{equation}\label{eq:G_deriv_approx_error}
\sup_{x \in [x_{0,M},x_{1,M}]} |G^{(j)}(x) - S_d^{(j)}(x)| = O(M^{-(p+1)+j}), \qquad \mbox{for}~j=0,1,2.
\end{equation}
Now, $S_d^{(1)}(x)$ is a spline of order $d-1$ on the same set of knots and hence can be expressed
as $\sum_{k=1}^M \beta_k^*\phi_{k,M}(x)$ for all $x \in [x_{0,M},x_{1,M}]$ if $\{\phi_{k,M}\}_{k=1}^M$
is the normalized spline basis of order $d-1$ on the same set of knots. Without loss of generality,
we assume that $x_{0,M} < x_{0,\delta} < x_{1,\delta} < x_{1,M}$. Then, by integration by parts, we have
\begin{eqnarray*}
\int_{x_{0,\delta}}^x \frac{g(u) - S_d^{(1)}(u)}{g(u)} dx &\hskip-.1in=& \hskip-.1in
\frac{1}{g(x)} (G(x) - S_d(x))   -
\frac{1}{g(x_{0,\delta})} (G(x_{0,\delta}) - S_d(x_{0,\delta})) \nonumber\\
&& ~~+ \int_{x_{0,\delta}}^x \frac{g'(u)}{(g(u))^2}
(G(u) - S_d(u))du.
\end{eqnarray*}
Since $(g(u))^{-1}$ and $g'(u)$ are bounded on $D$, (\ref{eq:integral_g_spline_approx}). As a by-product,
we also have from (\ref{eq:G_deriv_approx_error}) that $\sup_{x \in [x_{0,M},x_{1,M}]}| g^{(j)}(x) - g_{\bs\beta^*}^{(j)}(x)|
= O(M^{-p+j})$ for $j=0,1$. Moreover, it can be checked that if $p \geq 2$ then $\sup_{x \in [x_{0,M},x_{1,M}] \cap \Pi_M^c}
|S_d^{(3)}(x)|$ is bounded, where $\Pi_M$ constitute the knot sequence, which are the potential points of non-smoothness
for $S_d$.

\section{Proof of Theorem \ref{thm:optimal_rate}}\label{subsec:optimal_rate}

Since $\widehat{\bs\beta}$ is a
local minimizer of $\tilde L_\delta(\bs\beta)$, it satisfies $(\partial/\partial \bs\beta) \tilde L_\delta(\widehat{\bs\beta}) = 0$.
Thus, applying the Mean Value Theorem coordinatewise, we have, for $k=1,\ldots,M$,
\begin{eqnarray}\label{eq:L_tilde_gradient_expansion}
- \frac{\partial}{\partial \beta_k} \tilde L_\delta(\bs\beta^*) &=& \frac{\partial^2}{\partial \beta_k \partial \bs\beta^T}
\tilde L_\delta(\widetilde{\bs\beta}_k) (\widehat{\bs\beta} - \bs\beta^*)
\end{eqnarray}
for some $\widetilde{\bs\beta}_k$ such that $\parallel \widetilde{\bs\beta}_k - \bs\beta^*\parallel \leq
\parallel \widehat{\bs\beta} - \bs\beta^*\parallel$.
We write
\begin{eqnarray*}
- \frac{\partial}{\partial \bs\beta} \tilde L_\delta(\bs\beta^*)
&\hskip-.1in=& \hskip-.1in
2\sum_{j=1}^n \varepsilon_j \frac{\partial}{\partial \bs\beta} X(T_j;\bs\beta^*,\widehat x_0) \mathbf{1}_{[\delta,1-\delta]}(T_j)
\nonumber\\
&\hskip-.2in& \hskip-.2in
+ 2\sum_{j=1}^n (X_g(T_j) -  X(T_j;\bs\beta^*,\widehat x_0))\frac{\partial}{\partial \bs\beta} X(T_j;\bs\beta^*,\widehat x_0) \mathbf{1}_{[\delta,1-\delta]}(T_j) \nonumber\\
&=:& \tilde U_1 + \tilde U_2.
\end{eqnarray*}
On the other hand,
\begin{eqnarray*}
\frac{\partial^2}{\partial \bs\beta \partial \bs\beta^T}
\tilde L_\delta(\bs\beta) &=& 2\sum_{j=1}^n \frac{\partial}{\partial \bs\beta}
X(T_j;\bs\beta,\widehat x_0) \left(\frac{\partial}{\partial \bs\beta}
X(T_j;\bs\beta,\widehat x_0)\right)^T \mathbf{1}_{[\delta,1-\delta]}(T_j) \nonumber\\
&&
- 2\sum_{j=1}^n \varepsilon_j \frac{\partial^2}{\partial \bs\beta \partial \bs\beta^T} X(T_j;\bs\beta,\widehat x_0) \mathbf{1}_{[\delta,1-\delta]}(T_j)  \nonumber\\
&& - 2\sum_{j=1}^n (X_g(T_j) - X(T_j;\bs\beta,\widehat x_0))
\frac{\partial^2}{\partial \bs\beta \partial \bs\beta^T} X(T_j;\bs\beta,\widehat x_0)
\nonumber\\
&=:& \tilde S_1(\bs\beta) + \tilde S_2(\bs\beta) + \tilde S_3(\bs\beta).
\end{eqnarray*}
Then, we can express (\ref{eq:L_tilde_gradient_expansion})
in vectorial form as
\begin{eqnarray*}
- \frac{\partial}{\partial \bs\beta} \tilde L_\delta(\bs\beta^*)  &=&
\tilde S_1(\bs\beta^*)(\widehat{\bs\beta} - \bs\beta^*) + \sum_{k=1}^M \mathbf{e}_k\mathbf{e}_k^T
(\tilde S_1(\widetilde{\bs\beta}_k) - \tilde S_1(\bs\beta^*))(\widehat{\bs\beta} - \bs\beta^*)
\nonumber\\
&&
+ \sum_{k=1}^M \mathbf{e}_k\mathbf{e}_k^T
(\tilde S_2(\widetilde{\bs\beta}_k) + \tilde S_3(\widetilde{\bs\beta}_k))(\widehat{\bs\beta} - \bs\beta^*),
\end{eqnarray*}
where $\mathbf{e}_k$ denotes the vector in $\mathbb{R}^M$ with 1 in $k$-th coordinate and zero elsewhere. From
this, we get the expansion
\begin{eqnarray}\label{eq:beta_hat_expansion}
\widehat{\bs\beta} - \bs\beta^* &\hskip-.1in=& \hskip-.1in \left(\tilde S_1(\bs\beta^*)\right)^{-1}\tilde U_1 +
\left(\tilde S_1(\bs\beta^*)\right)^{-1} \tilde U_2 \nonumber\\
&\hskip-.1in& \hskip-.1in - \left(\tilde S_1(\bs\beta^*)\right)^{-1}\sum_{k=1}^M \mathbf{e}_k\mathbf{e}_k^T
(\tilde S_1(\widetilde{\bs\beta}_k) - \tilde S_1(\bs\beta^*))(\widehat{\bs\beta} - \bs\beta^*) \nonumber\\
&\hskip-.1in& \hskip-.1in
- \left(\tilde S_1(\bs\beta^*)\right)^{-1}\sum_{k=1}^M \mathbf{e}_k\mathbf{e}_k^T
(\tilde S_2(\widetilde{\bs\beta}_k) + \tilde S_3(\widetilde{\bs\beta}_k))(\widehat{\bs\beta} - \bs\beta^*). 
\end{eqnarray}
Let $S_l(\bs\beta)$ be the counterpart of $\tilde S_l(\bs\beta)$, $l=1,2,3$
once we replace the initial condition $\widehat x_0$ by $x_{0,\delta}$. It is then easily verified that
for some $C_{11} > 0$,
\begin{eqnarray}\label{eq:S_1_beta_diff_bound}
\parallel \frac{1}{n} \tilde S_1(\bs\beta^*) - \frac{1}{n} S_1(\bs\beta^*) \parallel 
&\leq&  C_{11}
M^{1/2} |\widehat x_0 - x_{0,\delta}| (1+ M^{1/2} |\widehat x_0 - x_{0,\delta}|) \nonumber\\
&=& O_P(M^{1/2}(\sigma_\varepsilon^2/n)^{(p+1)/(2p+3)}). 
\end{eqnarray}
This implies in particular, by Proposition \ref{prop:kappa} and Lemma \ref{lem:G_gamma_quad_bound} that
\begin{equation}\label{eq:S_1_beta_inverse_norm}
\parallel \left( \frac{1}{n} \tilde S_1(\bs\beta^*) \right)^{-1}  \parallel \leq C_{12} M^2
+  O_P(M^{1/2}(\sigma_\varepsilon^2/n)^{(p+1)/(2p+3)})
= O(M^2)
\end{equation}
for $M$ satisfying the condition (\ref{eq:M_alpha_optimal}).
Thus
\begin{equation*}
\mbox{Var}\left(\left(\tilde S_1(\bs\beta^*)\right)^{-1}\tilde U_1 ~|~\mathbf{T}, \widehat x_0\right) = \frac{\sigma_\varepsilon^2}{n}
\left(\frac{1}{n} \tilde S_1(\bs\beta^*)\right)^{-1},
\end{equation*}
and hence,
\begin{equation}\label{eq:phi_gradient_L_beta_x_variance}
\mbox{Trace} \left[\mbox{Var}\left(\left(\tilde S_1(\bs\beta^*)\right)^{-1}\tilde U_1 ~|~\mathbf{T}, \widehat x_0\right)\right]  = O_P\left(\frac{\sigma_\varepsilon^2 M^3}{n}\right)
\end{equation}
by (\ref{eq:S_1_beta_inverse_norm}).  On the other hand, from
\begin{equation*}
\sup_{t \in [\delta,1-\delta]} |X(t;\bs\beta^*,\widehat x_0) - X(t;\bs\beta^*)| = O(|\widehat x_0 - x_{0,\delta}|)
= O_P((\sigma_\varepsilon^2/n)^{(p+1)/(2p+3)})
\end{equation*}
and $\sup_{t \in [\delta, 1-\delta]}|X_g(t) - X(t;\bs\beta^*)| = O(M^{-(p+1)})$, and the form of $\tilde U_2$, we have,
with an application of Proposition \ref{prop:quadratic_solve} (stated below) that
\begin{eqnarray}\label{eq:phi_gradient_L_beta_x_bias}
\parallel (\tilde S_1(\bs\beta^*))^{-1} \tilde U_2 \parallel 
&\leq& 2 \parallel   \left(\frac{1}{n} \tilde S_1(\bs\beta^*)\right)^{-1}  \parallel^{1/2} \nonumber\\
&& ~~\cdot~
\sup_{t \in [\delta,1-\delta]}
(|X(t;\bs\beta^*,\widehat x_0) - X(t;\bs\beta^*)|  + |X_g(t) - X(t;\bs\beta^*)|) \nonumber\\
&=& O_P(M^{-p}) + O_P(M (\sigma_\varepsilon^2/n)^{(p+1)/(2p+3)}). 
\end{eqnarray}
Now, using arguments analogous to those used in the proof of Theorem \ref{thm:consistency},
and the bounds (\ref{eq:X_beta_bound})--(\ref{eq:X_diff_beta_beta_bound}), we can show that the maximum
of the norms of the matrices
$\left(\tilde S_1(\bs\beta^*)\right)^{-1}\sum_{k=1}^M \mathbf{e}_k\mathbf{e}_k^T
(\tilde S_1(\widetilde{\bs\beta}_k) - \tilde S_1(\bs\beta^*))$
and $\left(\tilde S_1(\bs\beta^*)\right)^{-1}\sum_{k=1}^M \mathbf{e}_k\mathbf{e}_k^T
(\tilde S_l(\widetilde{\bs\beta}_k)$,
for $l=2,3$, is $O_P(M^3 \alpha_n)$, which is $o_P(1)$ by the condition on $M$. From this, and
(\ref{eq:phi_gradient_L_beta_x_variance}) and (\ref{eq:phi_gradient_L_beta_x_bias}),
the result (\ref{eq:optimal_rate}) follows.

\begin{proposition}\label{prop:quadratic_solve}
Suppose that $B$ be an $p\times n$ matrix such that $BB^T$ is invertible. Let $y$ be an $n\times 1$ vector.
Then $\parallel (BB^T)^{-1} B y \parallel \leq (\parallel (BB^T)^{-1} \parallel)^{1/2} \parallel y \parallel$.
\end{proposition}
Proposition \ref{prop:quadratic_solve} follows immediately by using singular value decomposition of $B$.

\section{Rate of convergence of the two-stage estimator}\label{subsec:proofs_two_stage}

First, define
\begin{eqnarray*}
\mathbf{W} &=& \frac{1}{n}\sum_{j=1}^n \bs\phi(X(T_j))\bs\phi(X(T_j))^T\mathbf{1}_{[\delta,1-\delta]}(T_j)\\
\hat{\mathbf{W}} &=&
\frac{1}{n}\sum_{j=1}^n \bs\phi(\hat X(T_j))\bs\phi(\hat X(T_j))^T\mathbf{1}_{[\delta,1-\delta]}(T_j).
\end{eqnarray*}
Then, using the fact that $X'(t) = g(X(t))$ and $g_{\bs\beta^*}(x) = \bs\phi(x)^T \bs\beta^*$, we have
\begin{eqnarray}\label{eq:two_stage_beta_expansion}
\widetilde{\bs\beta} &\hskip-.1in=& \hskip-.1in\hat{\mathbf{W}}^{-1} \frac{1}{n} \sum_{j=1}^n \bs\phi(X(T_j))\bs\phi(X(T_j))^T \bs\beta^* \mathbf{1}_{[\delta,1-\delta]}(T_j) \nonumber\\
&\hskip-.1in& \hskip-.1in+ \hat{\mathbf{W}}^{-1}  \frac{1}{n} \sum_{j=1}^n (g(X(T_j)) - g_{\bs\beta^*}(X(T_j)))\bs\phi(X(T_j)) \mathbf{1}_{[\delta,1-\delta]}(T_j)\nonumber\\
&\hskip-.1in& \hskip-.1in + \hat{\mathbf{W}}^{-1}  \frac{1}{n} \sum_{j=1}^n (\hat X'(T_j) - X(T_j))\bs\phi(X(T_j)) \mathbf{1}_{[\delta,1-\delta]}(T_j)\nonumber\\
&\hskip-.1in& \hskip-.1in + \hat{\mathbf{W}}^{-1}\frac{1}{n} \sum_{j=1}^n \hat X'(T_j) (\bs\phi(\hat X(T_j)) - \bs\phi( X(T_j))) \mathbf{1}_{[\delta,1-\delta]}(T_j)\nonumber\\
&\hskip-.1in=& \hskip-.1in \bs\beta^* + \mathbf{W}^{-1}(\mathbf{W} - \hat{\mathbf{W}})\hat{\mathbf{W}}^{-1} \mathbf{W}\bs\beta^*  + R_1 + R_2 + R_3, 
\end{eqnarray}
where $R_1$, $R_2$ and $R_3$ are the expressions in the second, third and fourth lines after the first equality.

We check
that the following bounds hold with probability tending to 1 for any given sequence $M\to \infty$
as $n\to \infty$.
\begin{equation}\label{eq:W_condition_number}
\max\{\parallel \mathbf{W}\parallel, \parallel \mathbf{W}^{-1}\parallel\} = O(1).
\end{equation}
\begin{equation}\label{eq:phi_gram_bound}
\parallel \hat{\mathbf{W}} - \mathbf{W}\parallel = O(M^2 (\sigma_\varepsilon^2/n)^{(p+1)/(2p+3)}\sqrt{\log n}).
\end{equation}
\begin{eqnarray}
\parallel \mathbf{W}^{-1}\left(\frac{1}{n} \sum_{j=1}^n (g(X(T_j)) - g_{\bs\beta^*}(X(T_j))) \bs\phi(X(T_j)) \mathbf{1}_{[\delta,1-\delta]}(T_j)\right)\parallel  &=&  O(M^{-p}). \label{eq:g_diff_phi_bound}
\end{eqnarray}
\begin{eqnarray}
&& \parallel \mathbf{W}^{-1}\left(\frac{1}{n} \sum_{j=1}^n (\hat X'(T_j) - X'(T_j))\bs\phi(X(T_j))
\mathbf{1}_{[\delta,1-\delta]}(T_j)\right)\parallel \nonumber\\
&=&  O((\sigma_\varepsilon^2/n)^{p/(2p+3)}\sqrt{\log n}). \label{eq:X_prime_diff_phi_bound} \\
&&  \parallel \frac{1}{n} \sum_{j=1}^n \hat X'(T_j) (\bs\phi(\hat X(T_j))-\bs\phi(X(T_j))) \mathbf{1}_{[\delta,1-\delta]}(T_j)\parallel \nonumber\\
&=& O( M^{3/2} (\sigma_\varepsilon^2/n)^{(p+1)/(2p+3)}\sqrt{\log n}). \label{eq:X_prime_phi_diff_bound}
\end{eqnarray}
Combining these with (\ref{eq:two_stage_beta_expansion}) we obtain Proposition \ref{prop:two_stage_regression}.

\subsubsection*{Proof of (\ref{eq:W_condition_number})} First, write $\mathbf{W}$ as
$\bar{\mathbf{W}} + \Delta_W$, where
\begin{equation*}
\bar{\mathbf{W}} = \int_{\delta}^{1-\delta} \bs\phi(X(t))\bs\phi(X(t))^T f_T(t) dt\\
= \int_{X(\delta)}^{X(1-\delta)} \bs\phi(u)(\bs\phi(u))^T \frac{f_T(X^{-1}(u))}{g(u)}du.
\end{equation*}
Notice that for any $\mathbf{y} \in \mathbb{S}^{M-1}$, $\mathbf{y}^T  \bar{\mathbf{W}} \mathbf{y}$ lies in the interval
\begin{equation*}
\left(\int_{X(\delta)}^{X(1-\delta)} (\mathbf{y}^T \bs\phi(u))^2 du\right)
\left[\frac{\min_{s \in [0,1]}f_T(s)}{\max_{u \in  [X(0),X(1)]} g(u)}~,~
\frac{\max_{s \in [0,1]}f_T(s)}{\min_{u \in  [X(0),X(1)]}g(u)}\right]
\end{equation*}
from which it follows that $\max\{\parallel \bar{\mathbf{W}} \parallel, \parallel\bar{\mathbf{W}}^{-1}\parallel\}
= O(1)$ (uniformly in $M$) by the property of the B-spline basis (Schumaker, 2007). The result then follows
from the fact (derived along the line of Lemma \ref{lem:G_gamma_quad_bound})
that $\parallel \Delta_W\parallel \leq c(\eta) M\sqrt{\log n/n}
=o(1)$ with probability $1-n^{-\eta}$ for any given $\eta > 0$.

\subsubsection*{Proof of (\ref{eq:phi_gram_bound})} This follows from the observation that
\begin{eqnarray*}
&\hskip-.1in& \hskip-.1in\parallel \hat{\mathbf{W}} - \mathbf{W} \parallel \\
&\hskip-.1in\leq& \hskip-.1in \frac{1}{n} \sum_{j=1}^n (\parallel \bs\phi(\hat X(T_j))\parallel
+ \parallel \bs\phi(X(T_j))\parallel) \parallel \bs\phi(\hat X(T_j)) -  \bs\phi(X(T_j))\parallel\mathbf{1}_{[\delta,1-\delta](T_j)},
\end{eqnarray*}
and then using Mean Value Theorem, followed by
%$\sup_{x\in [x_{0,M},x_{1,M}]} \sum_{k=1}^M |\phi_{k,M}^{(j)}(x)|^2 = O(M^{1+2j})$ for $j=0,1,2$,
condition (iii) of {\bf A2}, and finally invoking (\ref{eq:rate_local_polynomial}), we get the result.

\subsubsection*{Proof of (\ref{eq:g_diff_phi_bound})} Here, if we denote the vector inside
$\parallel \cdot \parallel$ by $\bs\gamma$, then we have
$$
\mathbf{W}\bs\gamma = \frac{1}{n} \sum_{j=1}^n (g(X(T_j)) - g_{\bs\beta^*}(X(T_j))) \bs\phi(X(T_j))\mathbf{1}_{[\delta,1-\delta](T_j)}.
$$
Taking inner product with $\bs\gamma$, applying Cauchy-Schwarz inequality on the right and then using (\ref{eq:g_beta_estimates}), we have
$\bs\gamma^T \mathbf{W}\bs\gamma \leq c_9 M^{-p} \sqrt{\bs\gamma^T \mathbf{W}\bs\gamma}$
for some $c_9 > 0$. Hence, by (\ref{eq:W_condition_number}), we have the result.

\subsubsection*{Proof of (\ref{eq:X_prime_diff_phi_bound})} It is similar to that of (\ref{eq:g_diff_phi_bound}) and uses
(\ref{eq:rate_local_polynomial_deriv}) rather than (\ref{eq:g_beta_estimates}).

\subsubsection*{Proof of (\ref{eq:X_prime_phi_diff_bound})} It uses similar arguments as
in the proof of (\ref{eq:phi_gram_bound}).

\section{Sub-Gaussian random variables}\label{sec:Appendix_C}

We summarize a few facts about sub-Gaussian random variables.
The following is a restatement of Lemma 5.5 of Vershynin (2011).
\begin{lemma}\label{lem:subgaussian}
A random variable $\xi$ is sub-Gaussian, if any of the following equivalent
conditions hold.
\begin{itemize}
\item[(1)]
$\mathbb{E}(e^{\xi^2/K_1^2}) < \infty$ for some $0 < K_1 < \infty$
\item[(2)]
$(\mathbb{E}(|\xi|^q))^{1/q}\leq K_2 \sqrt{q}$ for all $q\geq 1$, for some $0 <
K_2 <\infty$.
\end{itemize}
If moreover, $\mathbb{E}(\xi) = 0$, then the following is equivalent to (1) and
(2).
\begin{itemize}
\item[(3)] $\mathbb{E}(e^{t\xi}) \leq e^{t^2 K_3^2}$ for all $t \in \mathbb{R}$, for some $0 < K_3 <
\infty$.
\end{itemize}
\end{lemma}

Define the \textit{sub-Gaussian norm} of a random variable $\xi$ to be
\begin{equation}
\parallel \xi\parallel_{\psi_2} := \sup q^{-1/2} (\mathbb{E}|\xi|^q)^{1/q}.
\end{equation}
Clearly, by Lemma \ref{lem:subgaussian}, $\xi$ is a sub-Gaussian random
variable if and only if $\parallel \xi\parallel_{\psi_2} < \infty$.

One of the useful characteristics of sub-Gaussianity is that it is preserved
under linear combinations. Specifically, we have the following result.
\begin{lemma}\label{lem:subgaussian_linear}
(Lemma 5.9 in Vershynin (2011)). Suppose that $X_1,\ldots,X_n$ are independent
sub-Gaussian random variables and $b_1,\ldots,b_n \in \mathbb{R}$ are nonrandom
quantities. Then $\sum_{i=1}^n b_i X_i$ is sub-Gaussian and
\begin{equation}\label{eq:subgaussian_norm_linear}
\parallel \sum_{i=1}^n b_i X_i \parallel_{\psi_2}^2 \leq C \sum_{i=1}^n b_i^2
\parallel X_i \parallel_{\psi_2}^2
\end{equation}
for some $C > 0$.
\end{lemma}
The result follows easily from the equivalent characterizations in Lemma
\ref{lem:subgaussian}, specifically, by using the moment generating function.
The following simple corollary is very useful for our applications.
\begin{corollary}\label{cor:subgaussian_linear_uniform}
Suppose that $X_1,\ldots,X_n$ are independent random variables with
$\max_{1\leq i \leq n} \parallel X_i \parallel_{\psi_2} \leq K < \infty$. Then
$\sum_{i=1}^n b_i X_i$ is sub-Gaussian and
\begin{equation}\label{eq:subgaussian_norm_linear_uniform}
\parallel \sum_{i=1}^n b_i X_i \parallel_{\psi_2}^2 \leq C K^2 (\sum_{i=1}^n
b_i^2)
\end{equation}
for some $C > 0$.
\end{corollary}

The following (Proposition 5.10 in Vershynin (2011)) is a version of
\textit{Hoeffding's inequality} for sub-Gaussian random variables.
\begin{lemma}\label{lem:subgaussian_Hoeffding}
Let $\xi_1,\ldots,\xi_n$ be independent random variables satisfying
$\mathbb{E}(\xi_i) = 0$, and let $K := \max_{1\leq i \leq n}
\parallel \xi_i \parallel_{\psi_2} < \infty$. Then for any $b_1,\ldots,b_n \in
\mathbb{R}$ we have
\begin{equation}\label{eq:subgaussian_Hoeffding}
\mathbb{P}\left(|\sum_{i=1}^n b_i \xi_i| > t\right) \leq e \exp\left(-\frac{c
t^2}{K^2 \sum_{i=1}^n b_i^2}\right), \qquad \mbox{for all}~t > 0,
\end{equation}
for some $c > 0$.
\end{lemma}

\bibliographystyle{plain}

\clearpage
\newpage

\setcounter{figure}{0}
\renewcommand{\thefigure}{S.\arabic{figure}}

\begin{figure}[th]
\begin{center}
\includegraphics[width=2.5in, height=2.5in, bb = 100 0 700 700]{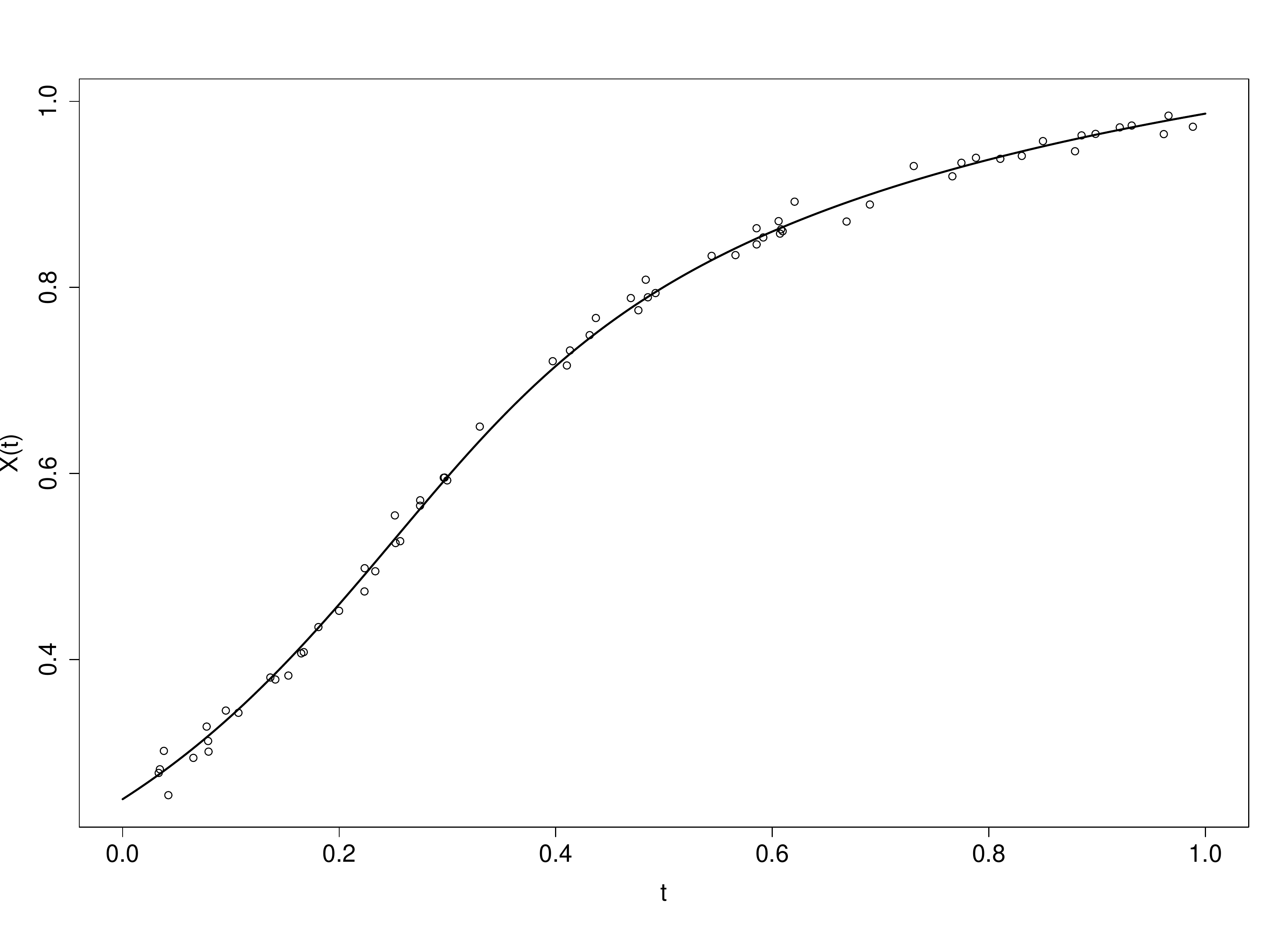}
\caption{Simulation: True trajectory and sample observations for one
replicate.
\label{fig:obs_trajectory}}
\end{center}
\end{figure}

\begin{figure}[th]
\begin{center}
\begin{tabular}{lllllllll}
\includegraphics[width=1.3in, height=1.5in, bb = 0 0 600
600]{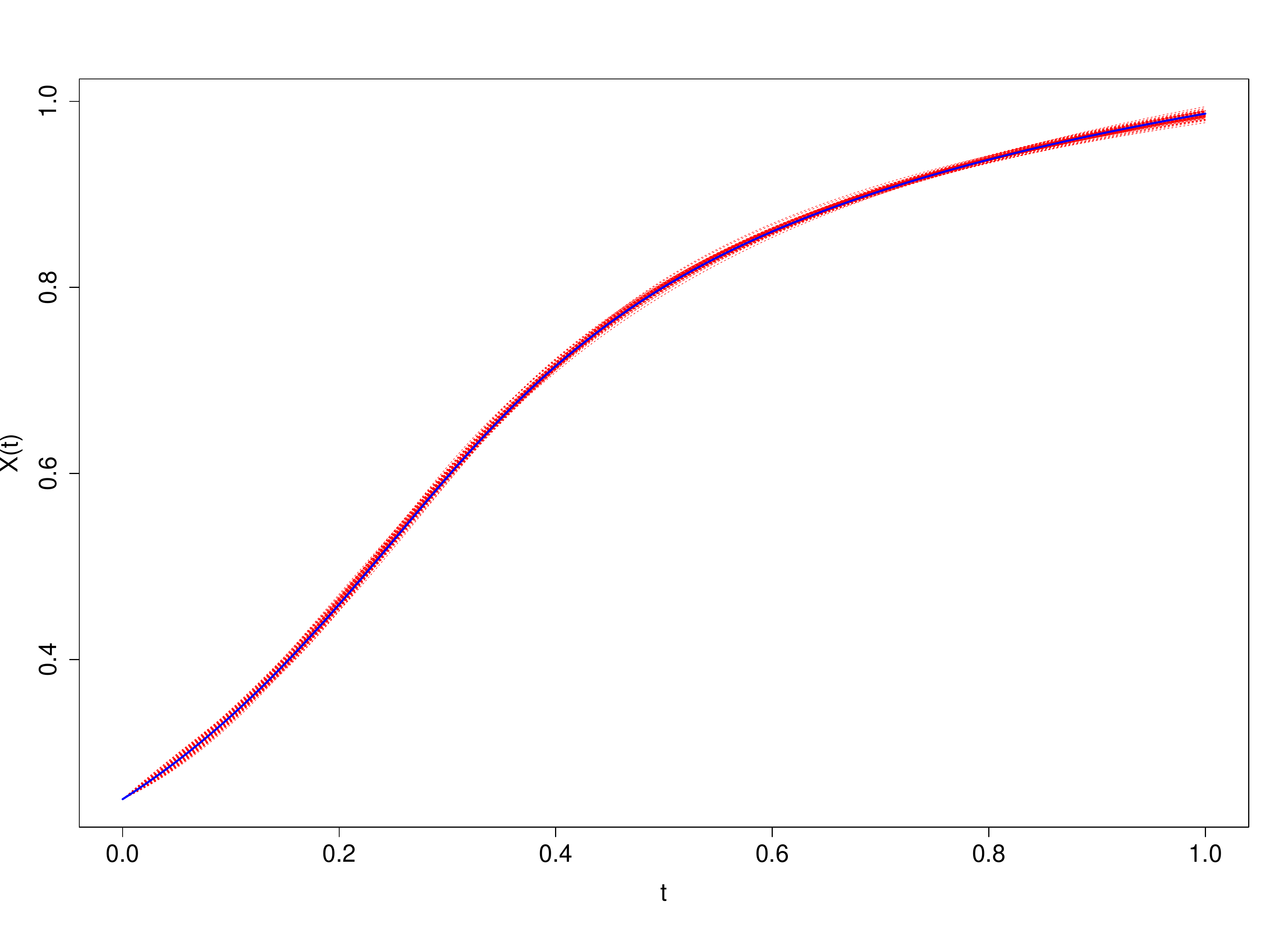} &&&&&&
\includegraphics[width=1.3in, height=1.5in, bb = 0 0 600
600]{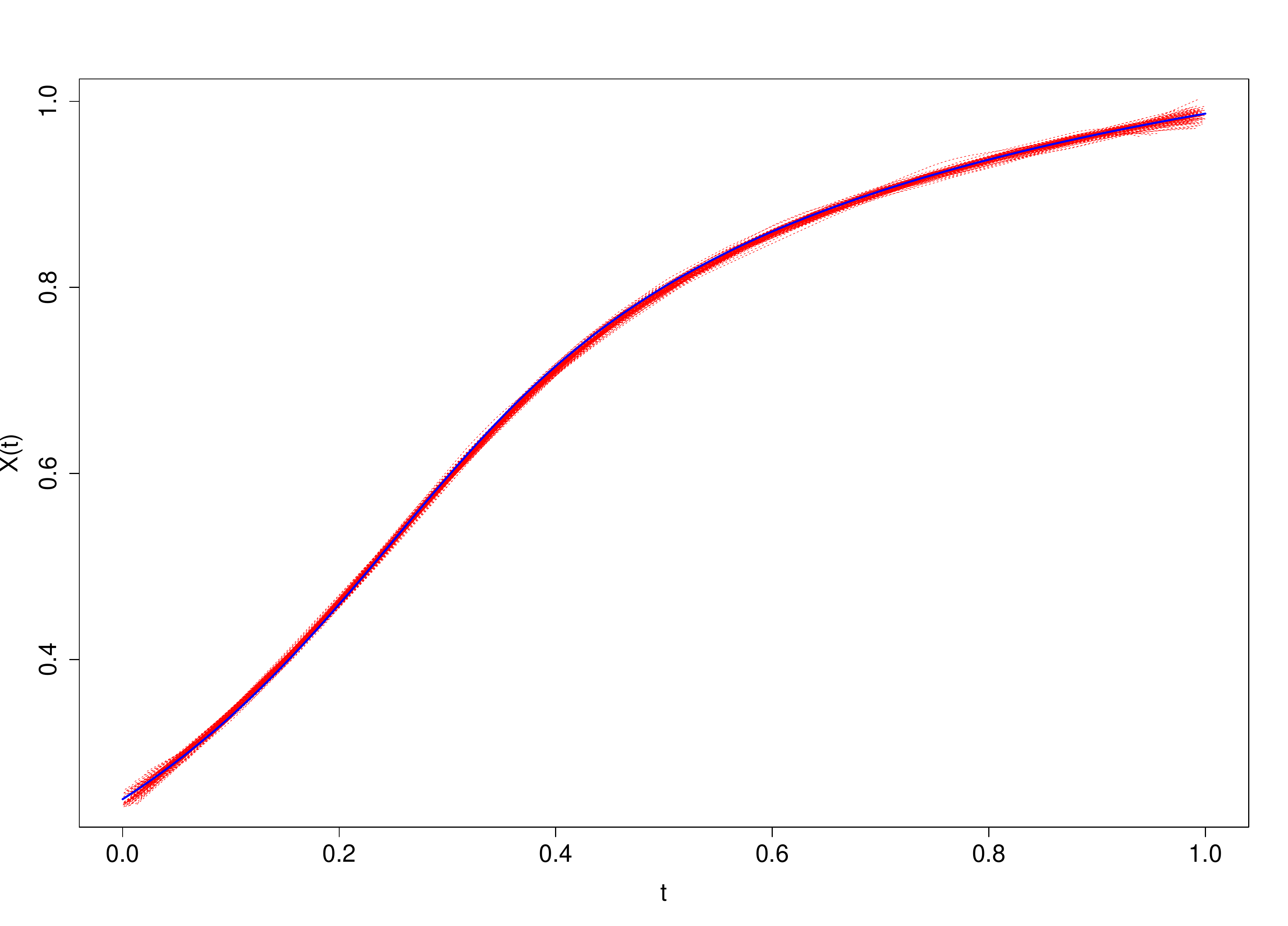}\\
\end{tabular}
\caption{Simulation: Estimated trajectory  $\widehat{X}(\cdot)$ (red curves) overlayed on the
true trajectory $X(\cdot)$ (blue curve). Left panel: proposed estimator; Right panel:
estimator from the 1st stage smoothing of the two-stage procedure.
\label{fig:simu_traj}}
\end{center}
\end{figure}

\begin{figure}[th]
\begin{center}
\begin{tabular}{lllllll}
\includegraphics[width=1.3in, height=1.5in, bb = 0 0 600
600]{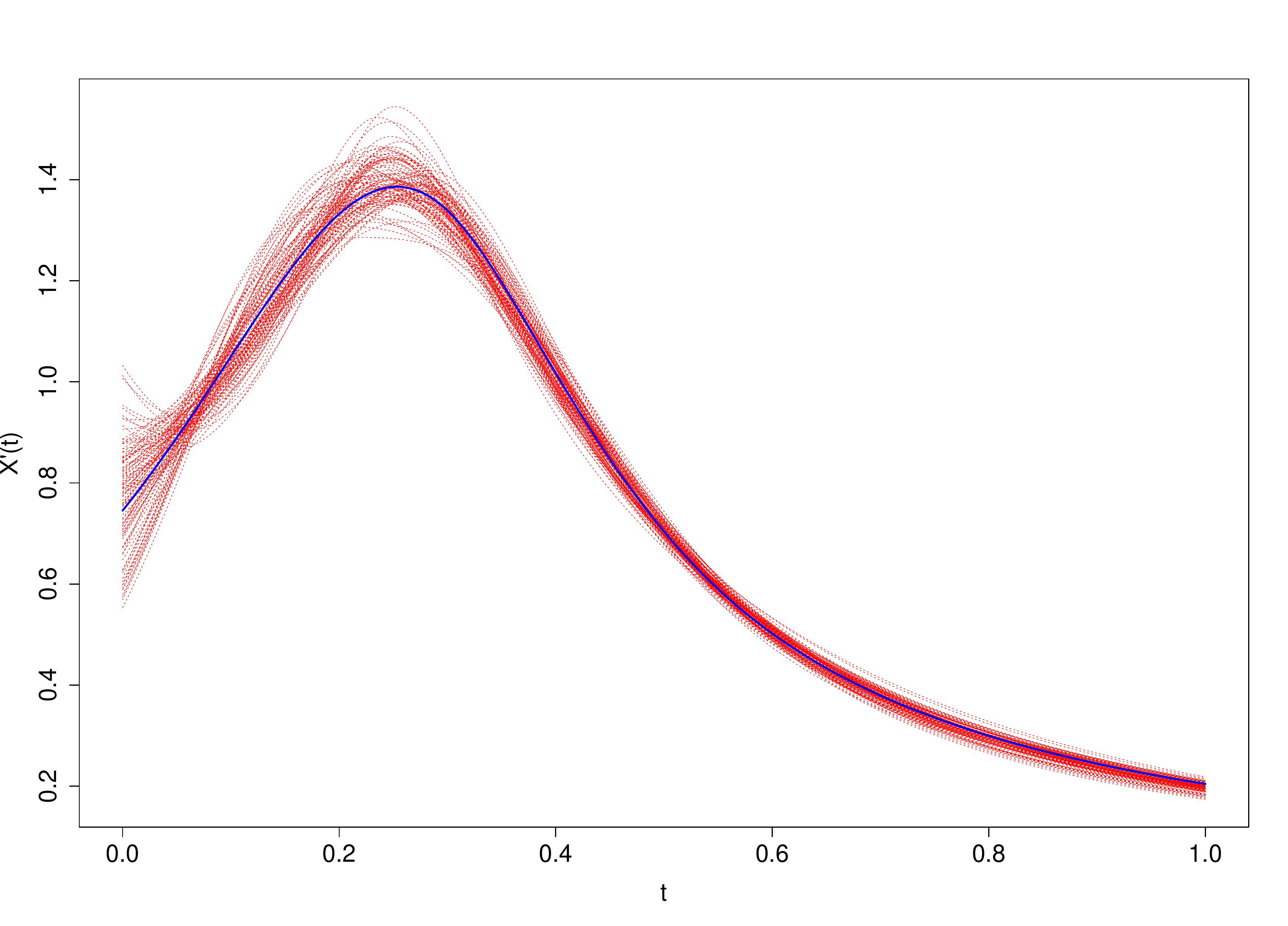} &&&&&
\includegraphics[width=1.3in, height=1.5in, bb = 0 0 600
600]{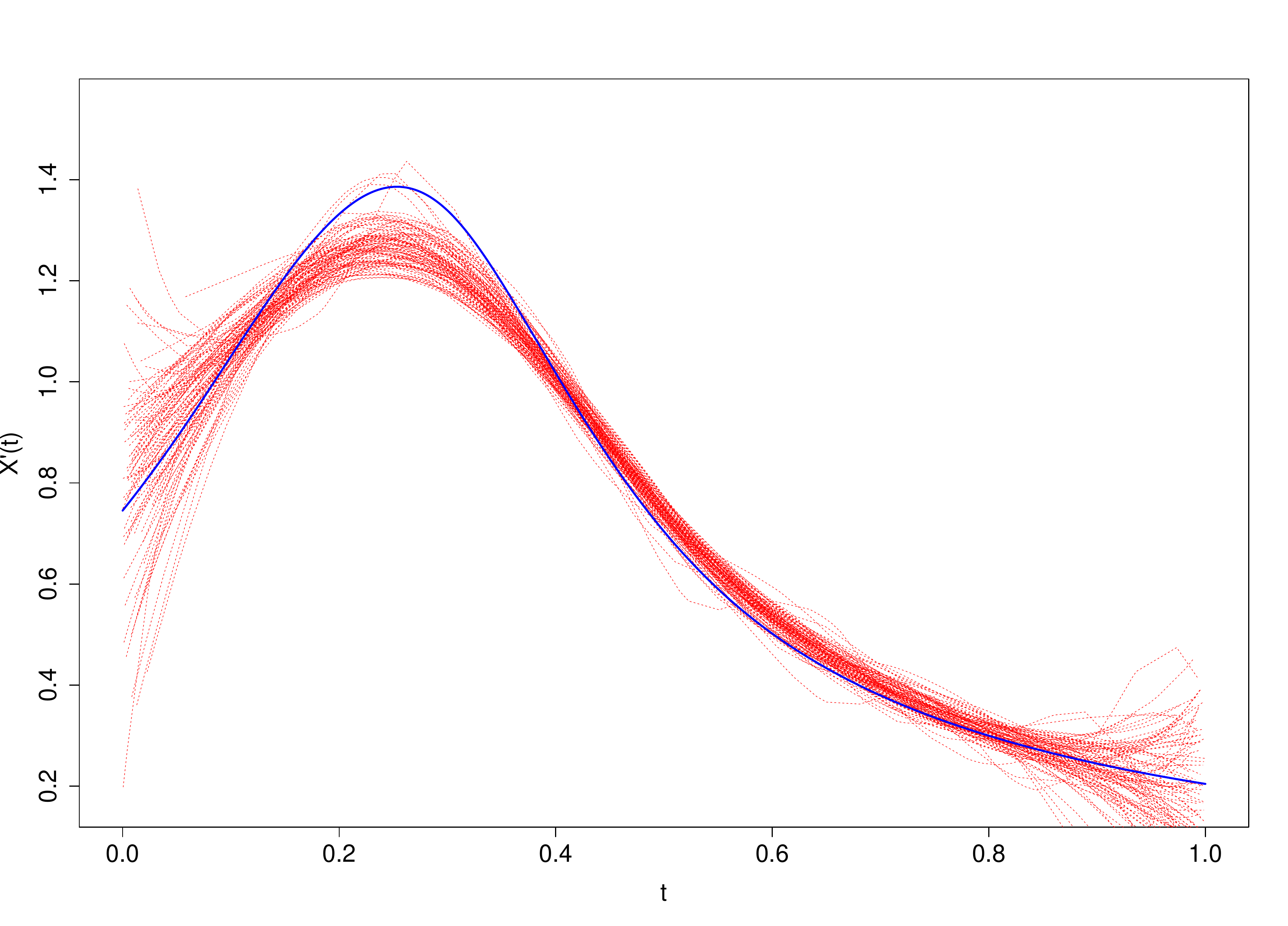}\\
\end{tabular}
\caption{Simulation: Estimated derivative of the trajectory $\widehat{X}^{\prime}(\cdot)$  (red curves) overlayed on the
true derivative  of  the trajectory $X^{\prime}(\cdot)$ (blue curve). Left panel: proposed estimator; Right panel:
estimator from the 1st stage smoothing of the two-stage procedure.
\label{fig:simu_traj_deriv}}
\end{center}
\end{figure}

\begin{figure}[th]
\begin{center}
\includegraphics[width=3.8in, height=4in, bb = 100 0 800 800]{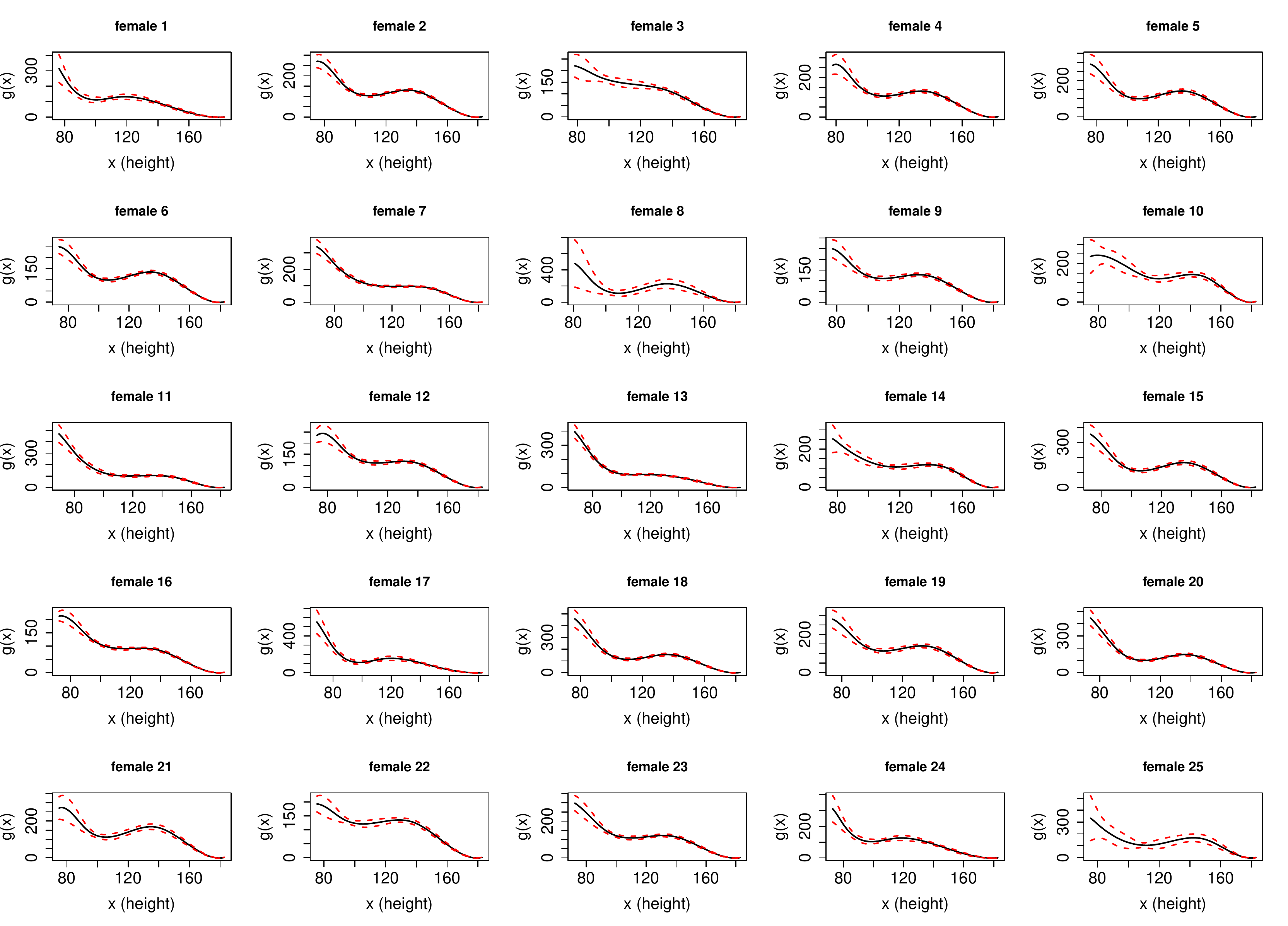}
\caption{Berkeley Growth Data: fitted gradient functions (black curve) for $25$
female subjects with two-standard-error bands (red broken lines).
\label{fig:fit_gradient_SE}}
\end{center}
\end{figure}

\begin{figure}[th]
\begin{center}
\includegraphics[width=3.8in, height=4in, bb = 100 0 800 800]{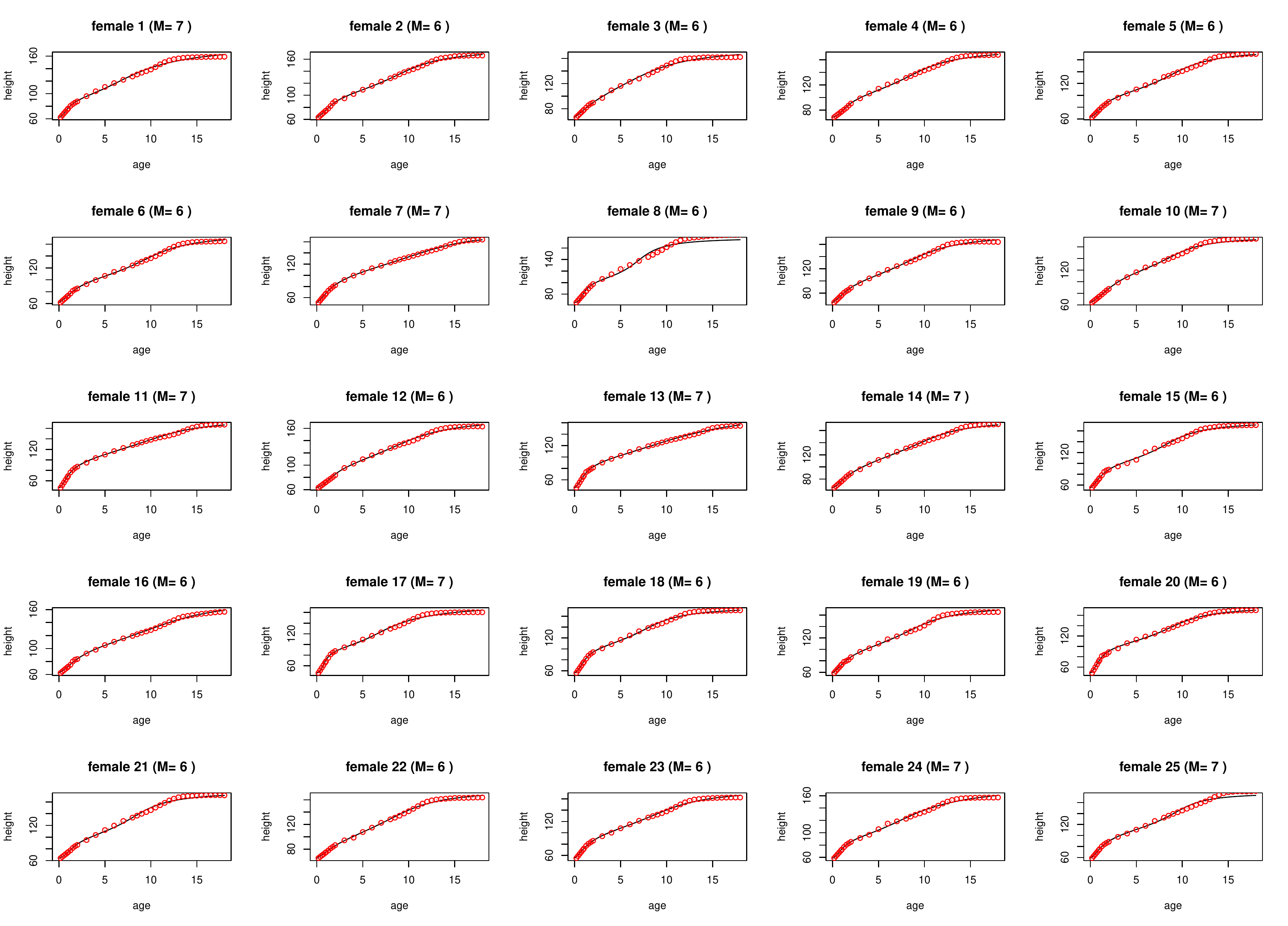}
\caption{Berkeley Growth Data: observed (red dots) and fitted (black curve)
growth trajectories for $25$ female subjects.
 \label{fig:fit_trajectory}}
\end{center}
\end{figure}

\end{document}